\documentclass[12pt,thmsa]{article}
\usepackage{graphicx}
\usepackage{amsmath}
\usepackage{amsfonts}
\usepackage{mathrsfs}
\usepackage{amssymb}

\newtheorem{Lemma}{Lemma}[section]{\bfseries}{\itshape}
\newtheorem{Pro}{Proposition}
{\bfseries}{\itshape}
\newtheorem{Remark}{Remark}
\newtheorem{Assumpt}{Assumption}

\newtheorem{Theorem}{Theorem}{\bfseries}{\itshape}

\newcommand{\R}{\mathbb R}
\newcommand{\Z}{\mathbb Z}
\newcommand{\N}{\mathbb N}
\newcommand{\p}{\mathbb P}
\newcommand{\E}{\mathbb E}

\begin{document}

\begin{center}
{\Large {\bf Wavelet-based estimation in a semiparametric regression model}}

{\Large \bigskip}

Emmanuel de Dieu NKOU and Guy Martial \ NKIET

\bigskip

Universit\'{e} des Sciences et Techniques de Masuku

BP 943 \ Franceville, Gabon

E-mail : emmanueldedieunkou@gmail.com,  gnkiet@hotmail.com.

\bigskip
\end{center}

\noindent\textbf{Abstract.}  In this paper, we introduce a wavelet-based method for estimating the EDR space in Li's  semiparametric regression model for achieving dimension reduction. This method is obtained by using linear wavelet estimators of the density and regression functions that are  involved in the covariance matrix of conditional expectation whose  spectral analysis gives the EDR directions. Then, consistency of the proposed estimators is proved. A simulation study that allow one to evaluate the performance of the proposal with comparison to existing methods is presented.  

\bigskip



\section{Introduction}\label{introduction}
Modeling by regression models allowing to establish links between a response variable and several explanatory variables is an approach that is both old and important in statistical analysis. In this perspective, many kinds of regression models have been introduced in the statistical literature, and the estimation problems related to these models have been intensively studied. For achieving both estimation and dimension reduction, a semiparametric regression model having the form 
\begin{equation}\label{pas1}
Y=F(\beta_1^{T}X,...,\beta_{N}^{T}X,\varepsilon), 
\end{equation}
have been introduced in \cite{Li1}. In this model, $Y$ is a real response random variable,  $X$ is a $d$-dimensional random vector containing the explanatory variables, $N$ is an integer of $\mathbb{N}^\ast$ such that $N<d$, $\beta_1,\cdots,\beta_N$ are unknown vectors in $\mathbb{R}^d$,  $\varepsilon$ is a real random variable that is independent  of $X$, and $F$ is an arbitrary unknown function.  It expresses the fact that the projection of $X$ onto the   subspace of $\mathbb{R}^d$  spanned by   $\beta_1,\cdots,\beta_{N}$, named the effective dimension reduction (EDR) space,  contains all information about   $Y$. Estimating $N$ and the EDR space  is then a crucial issue that have been tackled in several works (\cite{Aragon},\cite{Bura}, \cite{Duan1}, \cite{Ferre}, \cite{Hsing}, \cite{Li1}, \cite{Nkiet}, \cite{Nkou}, \cite{Saracco2}, \cite{Schott}, \cite{Velilla}, \cite{Zhu1}). For estimating the EDR space, it is enough to estimate the $\beta_j$'s which  are characterized, under some conditions, as  eigenvectors of the covariance matrix $\Lambda$    of $\mathbb{E}(X\vert Y)$.  For doing that, Li\cite{Li1} proposed a method,  called sliced inverse regression (SIR), based on slicing the range of $Y$. Although there exist alternative methods, this  method  stills the most popular  for dimension reduction. Based on the fact that the aforementioned  matrix is ​​expressed as a function of the density of $Y$ and regression functions,  Zhu and Fang \cite{Zhu1} proposed a nonparametric estimation procedure by using   kernel estimates of the preceding density and regression  functions.  But, as it is well known, there exist alternative nonparametric estimators for these functions. Among them, the wavelet-based estimators are known to have interesting properties and  have been successfully used in many fields
of Statistics (see, e.g. \cite{Vidakovic}).  However,  they never have been used for estimation in the model (\ref{pas1}). That is why,  we propose in this paper an estimation method for the EDR space related to this model, based on wavelet-based estimates of the density and regression functions involved in $\Lambda$. The rest of the paper is organized as follows.  In Section \ref{notations},  we construct an estimator  of $\Lambda$  based on wavelet method. Consistency of the resulting estimators are then given in   Section \ref{assumptions}.  Section \ref{simulations} is devoted to a simulation study made in order to evaluate the performance of the proposal with comparison to existing estimation methods.   The proofs of the main results are postponed in Section \ref{proofs}.
\medskip

\section{Wavelet-based estimation}\label{notations}
Putting $X=\left(X_1,\cdots,X_d\right)^T$ and letting  $f$ be the density of $Y$, we  suppose  that   for all $y\in\mathbb{R}$, we have  $f(y)>0$; then,    for any $j=1,\cdots,d$, we consider
\[
R_j(y)=\mathbb{E}(X_j\vert Y=y)=\frac{g_j(y)}{f(y)}
\,\,\textrm{ where }\,\,g_j(y)=\int_\mathbb{R} x f_{_{(X_j,Y)}}(x,y)dx,
\] 
\[
R(Y)=\left(R_1(Y),...,R_d(Y)\right)^T=\left(\mathbb{E}(X_1\vert Y),...,\mathbb{E}(X_d\vert Y)\right)^T=\mathbb{E}\left(X\vert Y\right),
\]
where $f_{_{(X_j,Y)}}$ denotes  the bivariate  density of the pair $(X_j,Y)$. The covariance matrix  
$
\Lambda = Cov\left(\mathbb{E}\left(X\vert Y\right)\right)
$
is of  a   great importance since the EDR space, which is to be estimated,    is obtained from its spectral analysis\cite{Li1}.  More precisely, $\beta_j$ is taken as an eigenvector of $\Lambda$ associated with the $j$-th largest eigenvalue $\lambda_j$. For estimating this matrix  Li \cite{Li1} used an approach based on slicing the range of $Y$ whereas Zhu and Fang\cite{Zhu1} introduced an estimator based on kernel estimates of the involved density and regression functions.  Here,  we propose an estimator of $\Lambda$ obtained from wavelets-based  estimates of   $f$ and $g_j$, $j=1,\cdots,d$. We assume that these functions belong to the space $\mathscr{L}^2(\R)$ of square integrable functions from $\mathbb{R}$ to itself.
Let $\varphi\in \mathscr{L}^2(\R)$ be a   father wavelet, and $\psi$ the associated mother wavelet, so that 
$\left\{\varphi_k(\cdot)=\varphi(\cdot-k),\, \psi_{\ell k}(\cdot)=2^{\ell /2}\psi\left(2^{\ell}(\cdot)-k\right)\,:\, k\in\Z,\,\ell\in\N\right\}$ is an orthonormal basis of $\mathscr{L}^2(\R)$ (ee, e.g., \cite{Daubechies1}, \cite{Hardle1}, \cite{Meyer2}).  Considering an i.i.d. sample  $\{(X^{(i)},Y_i)\}_{i=1,...n}$   of     $(X,Y )$,   and putting
$
X^{(i)}=\left(X_{i1},\cdots,X_{id}\right)^T,
$
we make use of the following estimators of $f$ and $g_j$:
\begin{equation*}\label{pas2}
\widehat{f}_n(y)= \sum_{k\in\Z}\widehat{\alpha}_{k}^{(n)}\varphi_{k}(y)+\sum_{\ell=0}^{j_n}\sum_{k\in\Z}\widehat{\gamma}_{\ell k}^{(n)}\psi_{\ell k}(y), 
\end{equation*}
\begin{equation*}
\widehat{g}_{j,n}(y)= \sum_{k\in\Z}\widehat{\delta}_{j,k}^{(n)}\varphi_{k}(y)+\sum_{\ell=0}^{j_n}\sum_{k\in\Z}\widehat{\eta}_{j,\ell k}^{(n)}\psi_{\ell k}(y), 
\end{equation*}
where 
\[
\widehat{\alpha}_{k}^{(n)} = \frac{1}{n}\sum_{i=1}^n \varphi_{k}\left(Y_i\right),\,\,\, \widehat{\gamma}_{jk}^{(n)}=\frac{1}{n}\sum_{i=1}^n\psi_{jk}\left(Y_i\right),
\]
\[
\widehat{\delta}_{j,k}^{(n)} = \frac{1}{n}\sum_{i=1}^n X_{ij}\varphi_{k}\left(Y_i\right),\,\, \, \widehat{\eta}_{j,\ell k}^{(n)}=\frac{1}{n}\sum_{i=1}^n X_{ij} \psi_{\ell k}\left(Y_i\right),
\]
and $(j_n)_{n\in\mathbb{N}}$ is an increasing  sequence of integers such that    $j_n\nearrow +\infty$ as $n \rightarrow +\infty$. As it was already done\cite{Zhu1}, in order to avoid  small values in the denominator, we  consider
\[
f_{b_n}(y)=\max\big(f(y),b_n\big)\,\,\textrm{ and }\,\, \widehat{f}_{b_n}(y)=\max\left(\widehat{f}_n(y),b_n\right),
\]
where  $\left(b_n\right)_{n\in\mathbb{N}^\ast}$ is a   sequence of positive real numbers  such that $\displaystyle \lim_{n\rightarrow +\infty}(b_n)=0$, and we estimate the ratio    
$R_{b_n,j}(y)= g_j(y)/f_{b_n}(y)$ by  
$
 \widehat{R}_{b_n,j}(y)= \widehat{g}_{j,n}(y)/\widehat{f}_{b_n}(y)$.
Then,  putting
$
\widehat{R}_{b_n}(y)=\left(\widehat{R}_{b_n,1}(y),...,\widehat{R}_{b_n,d}(y)\right)^T,
$
we take as estimator of   $\Lambda$ the random matrix: 
\[
\widehat{\Lambda}_n=\frac{1}{n}\sum_{i=1}^{n}\widehat{R}_{b_n}(Y_i)\widehat{R}_{b_n}(Y_i)^T.
\]
An estimate of the EDR space is obtained from the spectral analysis of this matrix. Indeed, if  $\widehat{\beta}_j$ is   as an eigenvector of $\widehat{\Lambda}_n$ associated with the $j$-th largest eigenvalue $\widehat{\lambda}_j^{(n)}$, we estimate the EDR space by the subspace of $\mathbb{R}^d$  spanned by   $\widehat{\beta}_1,\cdots, \widehat{\beta}_N$.
\begin{Remark}
It is well known that, under some conditions, the preceding wavelets estimators have linear foms given by
\[
\widehat{f}_n(y)=\frac{2^{j_n}}{n}\sum_{i=1}^n K\left(2^{j_n}y,\,2^{j_n}Y_i\right)\,\,\,\textrm{ and }\,\,\,
 \widehat{g}_{j,n}(y)=\frac{1}{n}\sum_{i=1}^n X_{ij}2^{j_n}K\left(2^{j_n}Y_i,2^{j_n}y\right),
\]
where
\begin{equation}\label{kernel}
 K\left(x,\,y\right) = \sum_{k\in\Z}\varphi(x-k)\varphi(y-k).
\end{equation}
This is the case when compactly supported wavelets, such as the Haar and Daubechies wavelets (see, e.g., \cite{Vidakovic}) for example, are used. In this case, the sum given in (\ref{kernel}) reduces to  a finite one (see \cite{Hardle1}).
\end{Remark}

\section{Assumptions and main results}\label{assumptions}

In this section, we give the used    assumptions, then  the main results  which give consistency of  $\widehat{\Lambda}_n$ and the $\widehat{\beta}_j$'s  are stated.  

\begin{Assumpt}\label{ass1}
The random variable $X$ is bounded, i.e.  there exists   $G>0$ such that
$\Vert X\Vert_d\leqslant G$, where $\Vert\cdot\Vert_d$ is the usual euclidean norm of $\mathbb{R}^d$.
\end{Assumpt}

\begin{Assumpt}\label{ass3}The $g_j$'s and 
 $f$ are $3$-times  differentiable and their  third derivatives  satisfy  the following condition: there exists a
neighborhood of the origin, say  $U$, and a constant $c>0$ such that, for any  $u\in U$, 
\[
\left|f^{(3)}\left(y+u\right)-f^{(3)}\left(u\right)\right|\leqslant c|u|\,\,\,\textrm{  and  }\,\,\,\left|g_j^{(3)}\left(y+u\right)-g_j^{(3)}\left(u\right)\right|\leqslant c|u|,
\]
for $j=1,\cdots,d$.
\end{Assumpt}

\begin{Assumpt}\label{ass4}
For any  $j\in\{1,\cdots,d\}$  and any   $u\in U$,
 $
\left|R_j(y+u) -R_j(y) \right|\leqslant c|u|.
$
\end{Assumpt}

\begin{Assumpt}\label{ass6}
The father wavelet  $\varphi$ and the mother wavelet   $\psi$ are bounded and  compactly  supported.
\end{Assumpt}

\begin{Assumpt}\label{ass7}The kernel $K$ given in (\ref{kernel}) satisfies the following properties:
\begin{enumerate} 
\item \label{b4} $\left|K(x,y)\right|\leqslant \Phi\left(x-y\right)$, where $\Phi:\R\longrightarrow \R_+$ is a bounded, compactly supported and symmetric function satisfying:
\[
\int u^2\Phi^2(u)\,du<+\infty\,\,\textrm{ and }\,\,\int |u|^{k}\Phi(u)dv <+\infty\,\,\textrm{ for }\,\,k\in\{0,1,4\}. 
\]
\item $\forall x \in \R$, $\forall k \in \left\{1,\,2,\,3\right\}$, 
$\int K(x,y)(y-x)^k dy=0$. 
 
 \end{enumerate}
\end{Assumpt}

\begin{Assumpt}\label{ass8}
 When $n$ is large enough $2^{-j_n}\sim n^{-c_1}$ and $b_n\sim n^{-c_2}$ where $c_1$ and $c_2$ are the positive numbers satisfying: $c_1>0$, $0<c_2<1/10$ and $1/8+c_2/4<c_1<1/4-c_2$.
\end{Assumpt}

\begin{Assumpt}\label{ass9}
\begin{enumerate}[R 1.]
\item For all  $j\in \left\{1,\cdots, d\right\}$, $\E\left[R_j^2(Y)\right]<+\infty$,
\item $\sqrt{n}\,\E\left[\left|R_k(Y)R_\ell(Y)\right|\mathbf{1}_{\{f(Y)\leqslant a_n \}}\right] = o(1)$  for $1\leqslant k,\ell \leqslant d$ and any sequence $(a_n)_{n\in\mathbb{N}}$ such that $a_b\sim b_n $ as $n\rightarrow +\infty$.

\end{enumerate}
\end{Assumpt}

\begin{Assumpt}\label{ass10}
The eigenvalues $\lambda_1, \cdots,\lambda_d$ of $\Lambda$ satisfy: $\lambda_1  > \cdots > \lambda_d > 0$.
\end{Assumpt}

\medskip

\begin{Remark}\label{rmq2}
Assumptions \ref{ass3} and \ref{ass4}     are of a type which is classical  in   nonparametric statistics  literature. They were used in particular in some papers on nonparametric estimation of the model which is tackled in this paper (see \cite{Nkou}, \cite{Zhu1}).    Several wavelet functions  satisfy  Assumption \ref{ass6} (see \cite{Daubechies1}, \cite{Hardle1}). That is the case, for instance,  for the  Daubechies wavelets and the Haar wavelets. The first condition in  Assumption \ref{ass7}-(1)  just is the condition  $H$ introduced in H$\ddot{\textrm{a}}$rdle et al.\cite{Hardle1},  the third condition is the condition  $H(k)$ and Assumption \ref{ass7}-(2) is a part of the condition $M(k)$.  Note that if Assumption \ref{ass6} holds, then for any $x\in\mathbb{R}$, $\int K(x,y)\,dy=1$, and there exists a non-zero constant $D$ independent of $x$ such that the inequality   $\int K^2(x,y)\,dy\leqslant D^2$ holds\cite{Gine3}. In Assumption \ref{ass8},  one can take   $b_n=\min(a,n^{-c_2})$, where $a$ is a fixed strictly positive number and sufficiently small. This yields a more accurate estimation of $f$; see  more details in  Remark 3.1 of Nkou and Nkiet\cite{Nkou}.  Assumption \ref{ass10}-(2)  is of a kind that has already been used   in the literature (see \cite{Zhu1}).
\end{Remark}

\bigskip

\noindent We now present the main results.
\bigskip

\begin{Theorem}\label{theo1}
Under the  assumptions \ref{ass1} to  \ref{ass9}, as $n\rightarrow +\infty$, we have
\[
\sqrt{n}\left(\widehat{\Lambda}_n-\Lambda\right)\stackrel{\mathscr{D}}{\rightarrow}\mathcal{H},
\]
where $\stackrel{\mathscr{D}}{\rightarrow}$ denotes convergence in distribution, $\mathcal{H}$ is a random variable having a normal distribution,  in the space $\mathscr{M}_d(\mathbb{R})$ of $d\times d$ matrices, such that, for any $A=\left(a_{k\ell}\right)\in \mathscr{M}_d(\mathbb{R})$, $A\neq 0$,  one has   $\textrm{tr}\left(A^T\mathcal{H}\right) \leadsto \mathcal{N}(0,\sigma^2_{A})$ with:  
\begin{equation}\label{siga}
\sigma^2_{A}=Var\bigg(\sum_{k=1}^d\sum_{\ell=1}^d\frac{a_{k\ell}}{2}\Big(X_{\ell}R_k\left(Y\right)+X_k R_{\ell}\left(Y\right)\Big)\bigg).
\end{equation}
\end{Theorem}

\bigskip

\noindent From Theorem \ref{theo1}, we can derive the asymptotic normality of the   eigenvectors. For $(j,k)\in\{1,\cdots,d\}^2$, we put  $\beta_j=\left(\beta_{j1},\cdots,\beta_{jd}\right)^T$ and  we consider the random variable
\begin{equation}\label{w}
\mathcal{W}_{jk}=\left(\sum_{\stackrel{r=1}{r\neq j}}^d\frac{\beta_{rk}}{\lambda_j-\lambda_r}\right)\,\sum_{p=1}^d\sum_{q=1}^d\frac{\beta_{jp}\beta_{jq}}{2}\Big(X_{q}R_p\left(Y\right)+X_p R_{q}\left(Y\right)\Big).
\end{equation}
Then, we have:
\bigskip

\begin{Theorem}\label{theo2}
Under the  assumptions \ref{ass1} to  \ref{ass10},  we have
$
\sqrt{n}\left(\widehat{\beta}_j-\beta_j\right)\stackrel{\mathscr{D}}{\rightarrow}\mathcal{N}(0,\Sigma_j)
$, as $n\rightarrow +\infty$,
where $\Sigma_j$ is the $d\times d$ covariance matrix of the random vector $ \mathcal{W}_{j}=\left(\mathcal{W}_{j1},\cdots,\mathcal{W}_{jd}\right)^T$.
\end{Theorem}

\medskip

\section{Simulation results}\label{simulations}
In this section, we present results of  simulations that was made in order to evaluate  the performance of the introduced wavelet-based  method and to compare it with some existing methods. We estimate the  EDR directions corresponding to the following models with dimension $d=5$:
\bigskip

\noindent\textbf{Model 1:} $Y = X_1+X_2+X_3+X_4+\varepsilon$;

\medskip

\noindent\textbf{Model 2:} $Y = X_1\left(X_1+X_2+1\right)+\varepsilon$;

\medskip

\noindent\textbf{Model 3:} $Y = X_1\,\left(0.5+\left(X_2+1.5\right)^{2}\right)^{-1}+\varepsilon$.

\bigskip

\noindent Model 1 corresponds to $N=1$ and $\beta_1=(1,1,1,1,0)^T$ whereas for Model 2 and Model 3 we have $N=2$, and the EDR directions are, respectively,  $\beta_1=(1,0,0,0,0)^T$, $\beta_2=(1,1,0,0,0)^T$ and $\beta_1=(1,0,0,0,0)^T$, $\beta_2=(0,1,0,0,0)^T$. Each data set was obtained as follows: $X=(X_1,X_2,X_3,X_4,X_5)^T$ is generated from a multivariate normal distribution  $\mathcal{N}(0,\textrm{\textbf{I}}_5)$, where $\textrm{\textbf{I}}_5$ is the $5\times 5$ identity matrix, $\varepsilon$ is generated from a standard normal distribution and $Y$ is  computed according to the above  models.
 We simulated  100 independent replications  of  samples of size $n=500$ generated as indicated above, over which means and standard deviations of the $\widehat{\beta}_j$'s were computed, together with the means   of  squared cosines  between  $\widehat{\beta}_j$ and  $\beta_j$,  $j=1,2$, given by\cite{Li1}:
$$
R^2\left(\widehat{\beta}_j\right)=\frac{\left(\widehat{\beta}^T_j   \beta_j\right)^2}{\widehat{\beta}^T_j \widehat{\beta}_j\cdot \beta^T_j  \beta_j}.
$$
Four methods were used for estimating the EDR directions: sliced inverse regression (denoted by SIR) with number of slices equal to $H=5$, the kernel method (denoted by Kernel) with quadratic kernel   $$K(x)=0,9375 \left(1-x^2 \right)^2\mathbf{1}_{[-1,\,1]}(x)$$ and bandwith $h_n= n^{-0.2}\simeq 0.2885$, the wavelet-based method with Haar wavelets (denoted by Wavelet (H)) and the wavelet-based method with Daubechies  wavelets of order $2$ (denoted by Wavelet (D)). The wavelet-based methods were taken with resolution level $j_n=0$ and $b_n=0.01$.
Table 1  presents the obtained   means  and standard deviations  (in parentheses) of $$\widehat{\beta}_1 = (\widehat{\beta}_{11}, \,\widehat{\beta}_{12},\, \widehat{\beta}_{13},\, \widehat{\beta}_{14}, \,\widehat{\beta}_{15})^T$$  for Model 1, whereas the obtained means and standard deviation of $R^2\left(\widehat{\beta}_j\right)$, $j=1,2$, for all three models, are given in Table 2.  Boxplots showing these later are also given in Figures 1 to 3. 
As  it is seen in Table 1, the four methods  have   good behaviours, but the wavelet-based  method with  Daubechies  wavelets of order 2  seems to yield more accurate estimates. This fact is confirmed in  Table 2 and Figures 1 to 3. Indeed, it is seen in Table 3 that wavelet (D) correspond to   a standard deviation  relatively reduced  compared with ththat of the   other methods. Figures show that   its estimate of $\beta_1$ is very good like that of the  other methods, but it is very powerful in   estimating  $\beta_2$.

\begin{table} 
\setlength{\tabcolsep}{0.08cm} 
\renewcommand{\arraystretch}{1} 
\begin{center}
\begin{tabular}{ccccccccccc}
\hline\hline \\
   Method    &  & $\widehat{\beta}_{11}$ &  & $\widehat{\beta}_{12}$  &  & $\widehat{\beta}_{13}$ &   &  $\widehat{\beta}_{14}$ &   &$\widehat{\beta}_{15}$ \\
\hline \hline                                                                      \\
SIR              &     & 0.9806     &   & 0.9940  &  &0.9947  &    &0.9893  &   & 0.0076  \\
               &     & (0.0508)      &   & (0.0528)  &  &(0.0431)  &    &(0.0471)  &   & (0.0512)  \\
              &     &       &   &    &  &   &    &   &   &   \\
Kernel              &     & 0.9737     &   & 0.9866  &  &0.9993  &    &0.9992  &   & 0.0032  \\
               &     & (0.0455)      &   & (0.0466)  &  &(0.0486)  &    &(0.0459)  &   & (0.0514)  \\
              &     &       &   &    &  &   &    &   &   &   \\
Wavelet  (H)             &     & 0.9815     &   & 0.9836  &  &0.9914  &    &1.0029  &   &- 0.0017  \\
               &     & (0.0567)      &   & (0.0465)  &  &(0.0423)  &    &(0.0472)  &   & (0.0544)  \\
              &     &       &   &    &  &   &    &   &   &   \\
Wavelet  (D)             &     & 0.9916     &   & 0.9893  &  &0.9988  &    &0.9973  &   &- 0.0031  \\
               &     & (0.0339)      &   & (0.0318)  &  &(0.0336)  &    &(0.0384)  &   & (0.0383)  \\
              &     &       &   &    &  &   &    &   &   &   \\
\hline \hline \\
\end{tabular}
\end{center}
\centering \caption{\label{table:tab4}Means and standard deviations  of $\widehat{\beta}_1$ over 100 replications for Model 1 with $n=500$.}
\end{table}

\begin{table}\label{tb2}
\setlength{\tabcolsep}{0.08cm} 
\renewcommand{\arraystretch}{1} 
\begin{center}
\begin{tabular}{ccccccccccccccc}
\hline\hline \\
     &  & Model 1 & &  &  &    & Model 2 &   &  & &   &    &  Model 3 &  \\
    &  &   & &  &    & & &  &   & &   &    &   &  \\
 \cline{7-9} \cline{12-15}\\
   Method    &  & $R^2\left(\widehat{\beta}_1\right)$  &  & &  & $R^2\left(\widehat{\beta}_1\right)$  &  & $R^2\left(\widehat{\beta}_2\right)$ & & &   &  $R^2\left(\widehat{\beta}_1\right)$ &   &$R^2\left(\widehat{\beta}_2\right)$ \\
\hline \hline                                                                      \\
SIR              &     & 0.9986    & &   &   & 0.9357  &  &0.7658  &    & &  &0.9512  &   & 0.7624  \\
               &     & (0.0010)    & &   &   & (0.0852)  &  &(0.1555)  &  &  &  &(0.0649)  &   & (0.2041)  \\
              &     &      & &  &   &    &  & &  &   &    &   &   &   \\
Kernel            &     & 0.9987 &   &   &   & 0.9726  &  &0.8914  &  &  &  &0.9756  &   & 0.8722  \\
               &     & (0.0010)    &   & &   & (0.0190)  &  &(0.1261)  &   &  & &(0.0202)  &   & (0.1441)  \\
              &     &       &   &   & &  &  &   &  & &   &   &   &   \\
Wavelet  (H)             &     & 0.9986  &   &  &   & 0.9572  &  &0.8679  &  &  &   &0.9716    &   & 0.9023  \\
               &     & (0.0009)    &  &  &   & (0.0413)  &  &(0.0979)  &  &  &  &(0.0298)  &   & (0.0943)  \\
              &     &     & &   &   &    &  &  &  &   &    &   &   &   \\
Wavelet (D)             &    &   0.9994  & &   &   & 0.9740  &  &0.9604  &   &  &  &0.9775  &   &0.9501  \\
               &     & (0.0006)  &    &  &   & (0.0237)  &  &(0.0362)  &  & &   &(0.0228)  &   & (0.0447)  \\
              &     &      &  &  &   &  &  &   &  &   &    &   &   &   \\
\hline\hline  \\
\end{tabular}
\end{center}
\centering \caption{\label{table:tab4}Means and standard deviations  of $R^2\left(\widehat{\beta}_{j}\right)$ , $j=1,2$, over 100 replications with $n=500$.}
\end{table}

\vspace{9cm}

\begin{figure}[!h]
\centering
\includegraphics[scale=0.7,bb=320 430 160 40]{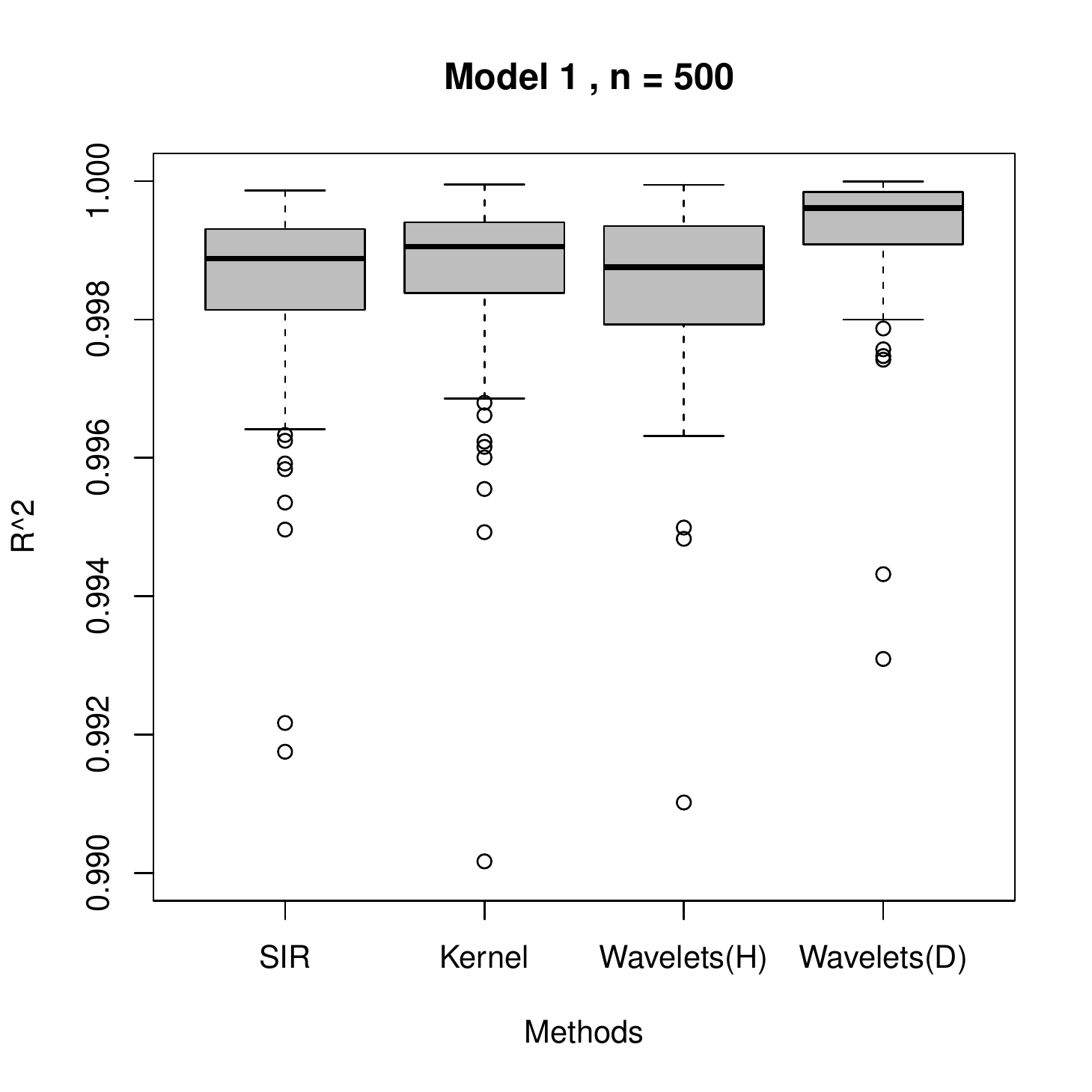}
\caption{Boxplots showing $R^2\left(\widehat{\beta}_1\right)$ for  Model 1,$ n = 500$.}
\label{boxplot_model1}
\end{figure}.


\vspace{6cm}

\begin{figure}[!h]
\centering
\leavevmode
\includegraphics[scale=0.5,bb=220 450 160 40]{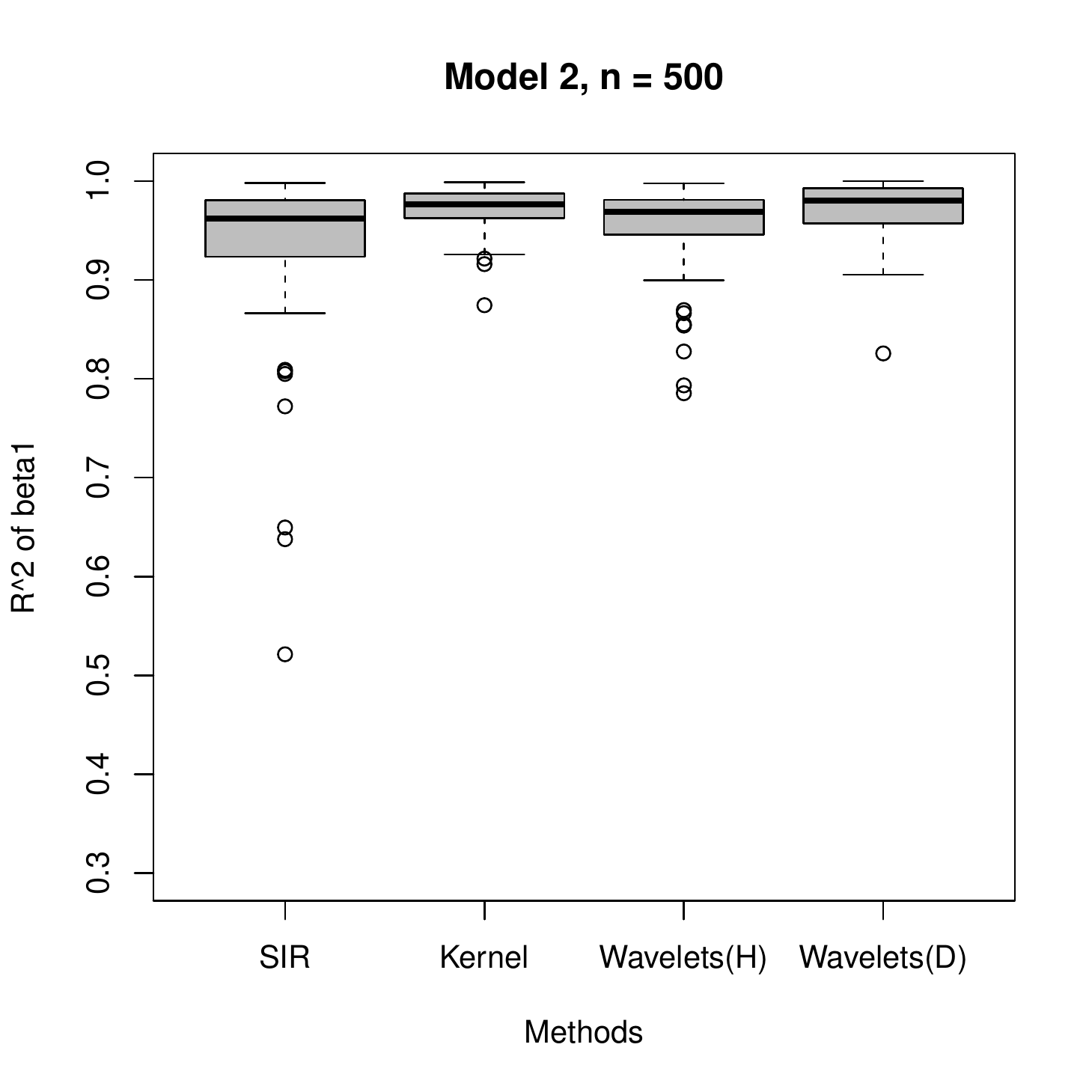}
\hspace{1cm}
\includegraphics[scale=0.5,bb=-180 450 160 40]{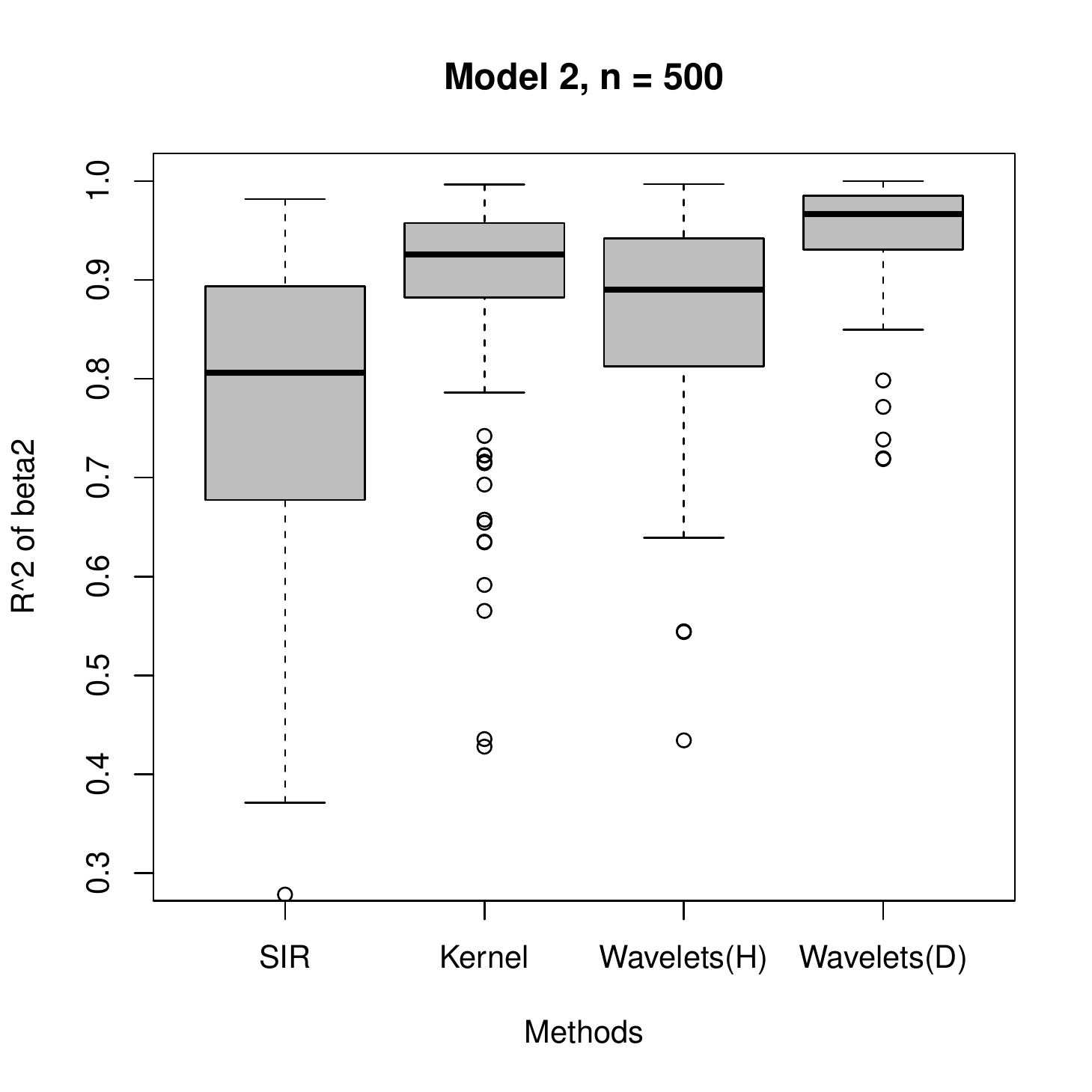}
\caption{Boxplots showing $R^2\left(\widehat{\beta}_j\right)$, $j=1,2$, for Model 2 with $n=500$ }\label{figure_model2}
\end{figure}.  

\vspace{10cm}

\begin{figure}[!h]
\centering
\leavevmode
\includegraphics[scale=0.5,bb=220 450 160 40]{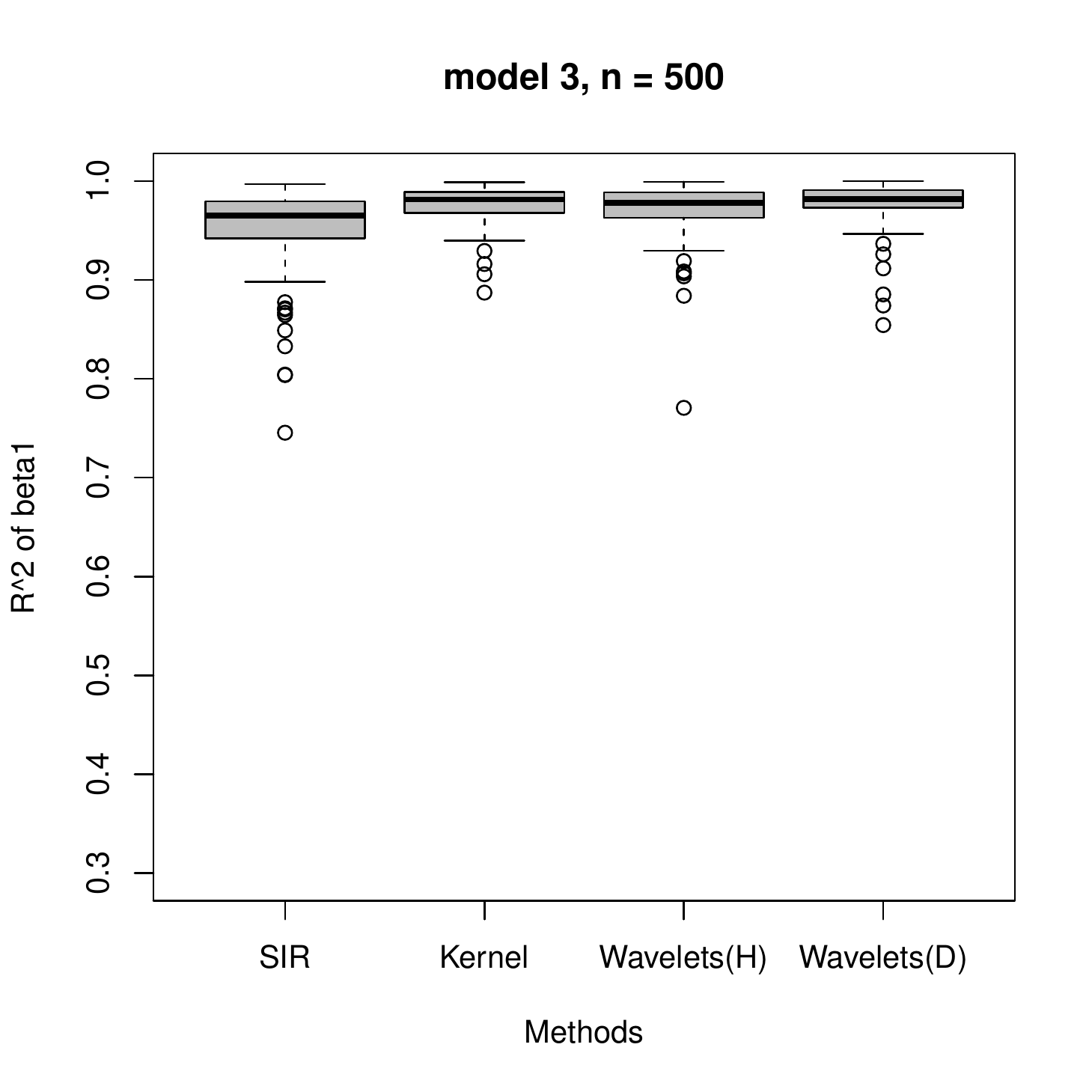}
\hspace{1cm}
\includegraphics[scale=0.5,bb=-180 450 160 40]{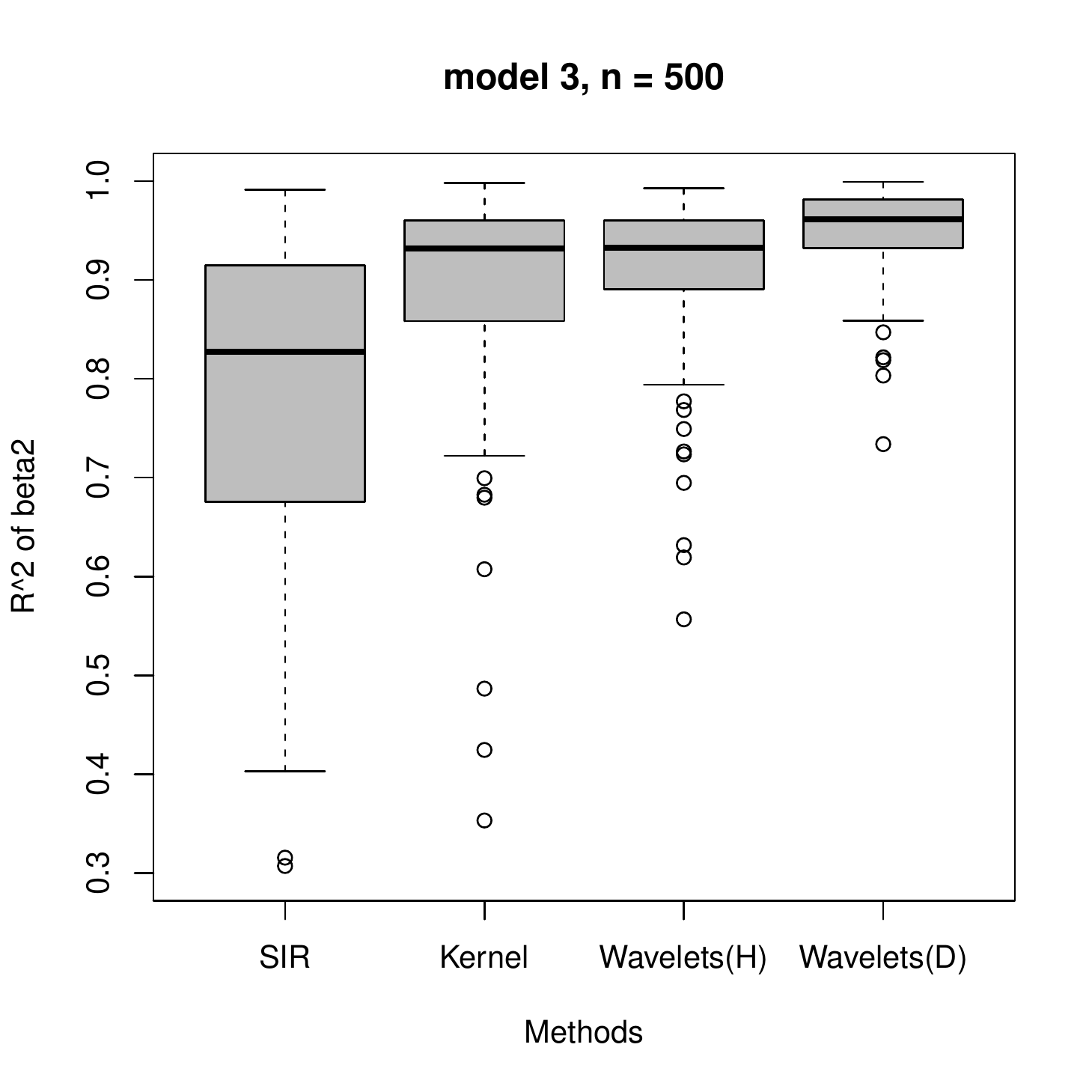}
\caption{Boxplots showing $R^2\left(\widehat{\beta}_j\right)$, $j=1,2$, for Model 3  with $n=500$. }\label{figure_model3}
\end{figure}.  

\vspace{7cm}

\section{Proofs}\label{proofs}

\subsection{Asymptotic properties of  $\widehat{f}_n$ and $\widehat{g}_{j,n}$}\label{asymptotics}

In this section, we give some results on asymptotic properties of $ \widehat{f}_n$ and $ \widehat{g}_{j,n}$ that are useful for proving the main result of the paper.  

\begin{Pro}\label{pro1}
Under the assumptions \ref{ass6} and \ref{ass7}, we have
\[
 \sup_{y\in\R}\left|\widehat{f}_n(y)-\E\left(\widehat{f}_n(y)\right)\right|=O\left(\left(\log n\right)^{\beta}\times 2^{j_n}n^{-1/2}\right) \, a.s.
\]
with $\beta>1/2$.
\end{Pro}
\noindent\textit{Proof.}
Consider the class  $
\mathcal{F}_n = \left\{\theta_y:u\mapsto \theta_y(u)=\frac{2^{j_n}}{n}K \left(2^{j_n}y,\,2^{j_n}u\right),\, y\in \R\right\}
$  of functions; each    $\theta_y$ is measurable since it is  a finite sum of measurables functions (Assumption \ref{ass6}).  By using the same  arguments than in  the proof of  Lemma 2 of   Gin\'e  and  Nickl\cite{Gine3}, it is easy to  check that    $\mathcal{F}_n$ is a $VC$- class of functions with respect to an  envelope $\theta$ such that $\left|\theta\right|\leqslant \frac{2^{j_n}}{n}\left\|\Phi\right\|_{\infty}$. 
Moreover,
\[
 \E\left(\left|\theta_y(\cdot)\right|\right)\leqslant \frac{2^{j_n}}{n}\left\|\Phi\right\|_{\infty}=:\mu_n \quad \mbox{ and }\quad \E\left(\left|\theta_y^2(\cdot)\right|\right)\leqslant \left(\frac{2^{j_n}}{n}\left\|\Phi\right\|_{\infty}\right)^2=\sigma_n^2.
\]
Then,  applying Talagrand's inequality (see Proposition 2.2 in Gin\'e and Guillou\cite{Gine2}; here $\mu_n=\sigma_n$), there
exist positive constants $K_1$ and $K_2$  such, that for all  
$t\geqslant K_1\left(\mu_n \log(A) + \sqrt{n}\sigma _n\sqrt{\log(A)}\right)$,
\begin{multline}
\p\Bigg\{\sup_{y\in\R}\left|\widehat{f}_n(y) - \E\left(\widehat{f}_n(y)\right)\right| > t\Bigg\}\nonumber\\
\leqslant  K_2 \exp\left\{-\frac{1}{K_2}\frac{t}{\mu_n}\log \left(1+\frac{t\mu_n}{K_2 \left(\sqrt{n}\sigma_n + \mu_n \sqrt{\log(A)}\right)^2}\right)\right\}\nonumber
=:v_n(t),
\end{multline}
where $A>0$. 
Taking $t:=t_n = \frac{2^{j_n}}{\sqrt{n}}\,\left(\log n\right)^{\beta}$, with $\beta >1/2$, we have   $\mu_n^{-1}\,t_n \rightarrow +\infty $  and $\sqrt{n}\mu_n\,t_n^{-1}\rightarrow 0$ as $n\rightarrow +\infty$. Thus, for $n$ large enough, $\mu_n^{-1}\,t_n-\sqrt{n}K_1\sqrt{\log (A)}\geqslant K_1\log (A)$, that is    $t_n\geqslant K_1\mu_n \left(\log (A) + \sqrt{n\,\log (A)}\right)$. Therefore, the preceding inequlity holds for $n$ large enough and, using $\log(1+u)\sim u$ and $(1+u)^2\sim 1$ as $u\rightarrow 0$, we get
\[
 v_n(t_n)\sim  K_2 \exp\left\{-\frac{\left(\log n\right)^{2\beta}}{ K_2^2 \left\|\Phi\right\|_{\infty}^2}\right\},
\]
from wich we deduce that $\sum_{n=1}^{+\infty}v_n(t_n)<+\infty$. Hence,
\begin{eqnarray*}\label{talagrand3}
\sum_{n=1}^{+\infty}\p\left\{\sup_{y\in\R}\left|\widehat{f}_n(y) - \E\left(\widehat{f}_n(y)\right)\right|  > t_n\right\}<+\infty,
\end{eqnarray*}
and   the required result  is obtained by using   Borel Cantelli's lemma.

\hfill $\square$

\begin{Pro}\label{pro2}
Under the assumptions \ref{ass3}, \ref{ass6} and \ref{ass7}, we have
\[
 \sup_{y\in\R}\left|\widehat{f}_n(y)-f(y)\right|=O\left(\left(\log n\right)^{\beta}\times 2^{j_n}n^{-1/2}\,+\,2^{-4j_n}\right)\mbox{ a.s. }
\] 
\end{Pro}
\noindent\textit{Proof.}
We can see that
\[
\E\left(\widehat{f}_n(y)\right)-f(y)=\int K\left(2^{j_n}y,2^{j_n}y+u\right)\left(f(y+u2^{-j_n})- f(y)\right)\,du,
\]
and by a Taylor expansion, we have: 
\begin{multline*}\label{b105}
\displaystyle
f(y+u2^{-j_n})-f(y)=\sum_{k=1}^2 \frac{2^{-kj_n}f^{(k)}(y)}{k!}u^k \,+\,\frac{2^{-3j_n}u^3}{2}\int_0^1(1-v)^{2}f^{(3)}\left(y+vu2^{-j_n}\right)dv.
\end{multline*} 
Then, using Assumption \ref{ass7}-(2) we get
\begin{multline}
\E\left(\widehat{f}_n(y)\right)-f(y) \nonumber\\
=\frac{2^{-3j_n}}{2}\int \int_0^1u^3(1-v)^{2}K\left(2^{j_n}y,2^{j_n}y+u\right)\bigg(f^{(3)}\left(y+vu2^{-j_n}\right) -f^{(3)}(y)\bigg) dv\,du,\nonumber
\end{multline}
and, therefore, 
\begin{eqnarray*}
&&\left|\E\left(\widehat{f}_n(y)\right)-f(y)\right| \nonumber\\
&\leqslant & \frac{2^{-3j_n}}{2}\int \int_0^1|u|^3(1-v)^{2}|K\left(2^{j_n}y,2^{j_n}y+u\right)| \left|f^{(3)}\left(y+vu2^{-j_n}\right) -f^{(3)}(y)\right| dv\,du\\
                                                  &\leqslant&  \frac{2^{-3j_n} c }{2}\int \int_0^1 |u|^{4}\Phi(u)  v2^{-j_n}(1-v)^{2}\,dv\,du\nonumber\\
				&\leqslant&  \frac{2^{-4j_n} c }{2}\int |u|^{4}\Phi(u) \, du.\nonumber
\end{eqnarray*}
We deduce that 
$\,\, \displaystyle
\sup_{y\in\R} \left|\E\left(\widehat{f}_n(y)\right)-f(y)\right|=O\left(2^{-4j_n}\right)
$. Combining this result to that of Proposition \ref{pro1} yields the required result.
\hfill $\square$ 

\begin{Pro}\label{pro3}
Under the assumptions \ref{ass1}, \ref{ass6} and \ref{ass7}, we have for $\beta>1/2$
\[
  \sup_{y\in\R}\left|\widehat{g}_{j,n}(y)-g_j(y)\right|=O\left(2^{j_n}n^{-1/2}\left(\log n\right)^{\beta} \, +\, 2^{-4j_n}\right) \, a.s..
\]
\end{Pro}
\noindent\textit{Proof.}  A similar reasoning than that of Proposition \ref{pro1} based on the class $\mathcal{G}_n=\left\{\theta_y:\left[-G;\, G\right]\times\R\longrightarrow \R \, / \, \theta_y(u,v)=2^{j_n}\frac{u}{n}K_{\varphi}\left(2^{j_n}y,2^{j_n}v\right), y\in\R\right\}$, where $G$ is the constant introduced in Assumption \ref{ass1},  gives 
\[
\sup_y\left|\widehat{g}_{j,n}\left(y\right)-\E\bigg(\widehat{g}_{j,n}\left(y\right)\bigg)\right|=O\left(2^{j_n}n^{-1/2}\left(\log n\right)^{\beta}\right) \,\,a.s..
\]
A similar reasoning than that of Proposition \ref{pro2} allows to obtain
$
\sup_{y\in\R} \left|\E\left(\widehat{g}_{j,n}(y)\right)\, -\,g_j(y)\right|=O\left(2^{-4j_n}\right)
$.
\hfill $\square$
\subsection{Preliminary lemmas}

\begin{Lemma}\label{lemma41} 
Under assumptions \ref{ass3}, \ref{ass6} and \ref{ass7},   we have 
\begin{equation}\label{b108}
\left|\int K \left(2^{j_n}x,2^{j_n}y\right)\left(g_{\ell}(x)-g_{\ell}(y)\right)\,dx \right|\, \leqslant \,C_1\,2^{-5j_n},
\end{equation}
where $C_1>0$.
\end{Lemma} 
\noindent\textit{Proof.}
From Taylor expansion of $g_{\ell}$ and   Assumption \ref{ass7}-(2) we get
\begin{multline*}
\int K\left(2^{j_n}x,2^{j_n}y\right)\left(g_\ell(x)-g_\ell(y)\right)dx\\ =  \frac{1}{2}\int\int_0^1 (x-y)^3 K\left(2^{j_n}x,2^{j_n}y\right)(1-v)^{2}\left(g_\ell^{(3)}(y+v(x-y)) -  g_\ell^{(3)}(y) \right)\,dv\,dx
\end{multline*}
and, therefore,
\begin{eqnarray*}
\left|\int K\left(2^{j_n}x,2^{j_n}y\right)\left(g_\ell(x)-g_\ell(y)\right)dx\right|&\leqslant & \frac{c}{2}\int |x-y|^{4} \Phi(2^{j_n}(x-y))dx\\
&=&\frac{2^{-5j_n}c}{2}\int |u|^{4}\Phi(u)\,\,du. 
\end{eqnarray*}
\hfill $\square$

\bigskip

\noindent Define, for $1\leqslant k, \ell\leqslant d$,
\begin{eqnarray*}
 U_{n,k,\ell}^{(1)} &=& \frac{1}{\sqrt{n}}\sum_{i=1}^n\Big\{g_k(Y_i)\left(\widehat{g}_{\ell , n}(Y_i)-g_{\ell}(Y_i)\right) + g_{\ell}(Y_i)\left(\widehat{g}_{k , n}(Y_i)-g_k(Y_i)\right)\Big\} \frac{\widehat{f}_{b_n}^2(Y_i)-f_{b_n}^2(Y_i)}{\widehat{f}_{b_n}^2(Y_i) f_{b_n}^2(Y_i)},\\
U_{n,k,\ell}^{(2)} &=& \frac{1}{\sqrt{n}}\sum_{i=1}^n \frac{\left(\widehat{g}_{k , n}(Y_i)-g_k(Y_i)\right)\left(\widehat{g}_{\ell , n}(Y_i)-g_{\ell}(Y_i)\right)}{\widehat{f}_{b_n}^2(Y_i)},\\
U_{n,k,\ell}^{(3)} &=& \frac{1}{\sqrt{n}}\sum_{i=1}^n R_{b_n,k}(Y_i)R_{b_n,\ell}(Y_i)\frac{\left(\widehat{f}_{b_n}^2(Y_i)-f_{b_n}^2(Y_i)\right)^2}{\widehat{f}_{b_n}^2(Y_i)f_{b_n}^2(Y_i)},\\
U_{n,k,\ell}^{(4)} &=& \frac{1}{\sqrt{n}}\sum_{i=1}^n\left(\widehat{f}_{b_n}(Y_i)-f_{b_n}(Y_i)\right)^2\frac{R_{b_n,k}(Y_i)R_{b_n,\ell}(Y_i)}{f_{b_n}^2(Y_i)}.
\end{eqnarray*}

\begin{Lemma}\label{cl14}
Under the assumptions \ref{ass1}, \ref{ass3}, \ref{ass6}, \ref{ass7} and \ref{ass8}, we have
$\left|U_n^{(q)}\right| = o_p(1)$ for any 
$q\in\{1,\cdots,4\}$.
\end{Lemma}
\noindent\textit{Proof.}
Using  Proposition \ref{pro3} and Assumption \ref{ass8}, we get
\begin{multline*}\label{b61}
\left|g_k(Y_i)\left(\widehat{g}_{\ell , n}(Y_i)-g_{\ell}(Y_i)\right) + g_{\ell}(Y_i)\left(\widehat{g}_{k , n}(Y_i)-g_k(Y_i)\right) \right| \\
    \leq A_1 \bigg(\left|g_k(Y_i)\right|+\left|g_{\ell}(Y_i)\right|\bigg)\, \tau_n ,
\end{multline*}    
where $A_1>0$ and  $\tau_n = n^{c_1-1/2}\left(\log n\right)^{\beta} \, +\, 2^{-4j_n}$. Then, using the inequalities  
$\left|g_k(Y_i)\right|+\left|g_{\ell}(Y_i)\right|\leqslant f_{b_n}(Y_i)\left(\left|R_k(Y_i)\right|+\left|R_{\ell}(Y_i)\right|\right)$,  
$\left|\widehat{f}_{b_n}(Y)-f_{b_n}(Y)\right| \leqslant \left|\widehat{f}_n(Y)-f(Y)\right|$,  $\widehat{f}_{b_n}(Y_i)\geq b_n$, $f_{b_n}(Y_i)\geq b_n$, together with Proposition \ref{pro2},  we obtain
\begin{eqnarray*}
& &\left|U_{n,k,\ell}^{(1)}\right|\\
& \leqslant & A_1\tau_n\,\frac{\sqrt{n}}{n}\sum_{i=1}^n\left[ \bigg(\left|R_k(Y_i)\right|+\left|R_{\ell}(Y_i)\right|\bigg)  \left|\frac{\left(\widehat{f}_{b_n}(Y_i)-f_{b_n}(Y_i)\right)^2}{\widehat{f}_{b_n}^2(Y_i)f_{b_n}(Y_i)} + 2\frac{\widehat{f}_{b_n}(Y_i) - f_{b_n}(Y_i)}{\widehat{f}_{b_n}^2(Y_i)}\right|\right]\\
                      & \leqslant &A_1A_2\, \tau_n^2\,\sqrt{n}\,b_n^{-2}\,\left(A_2b_n^{-1}\tau_n + 2\, \right) \frac{1}{n}\sum_{i=1}^n\bigg(\left|R_k(Y_i)\right|+\left|R_{\ell}(Y_i)\right|\bigg),
\end{eqnarray*}
where $A_2>0$. This implies that $
\left|U_{n,k,\ell}^{(1)}\right| = o_p(1)$,  due to the law of large numbers and the fact that:    $\lim_{n\rightarrow +\infty} \tau_n^2\,\sqrt{n}\,b_n^{-2}\,\left(b_n^{-1}\tau_n + 2\, \right) =0$. By a similar reasoning, we get
\[
\left|U_{n,k,\ell}^{(2)}\right|\leqslant n^{-1/2}b_n^{-2}\sum_{i=1}^n\left|\widehat{g}_{k , n}(Y_i)-g_k(Y_i)\right|\left|\widehat{g}_{\ell , n}(Y_i)-g_{\ell}(Y_i)\right|
\leq A_1^2n^{1/2}\,b_n^{-2}\,\tau_n^2;
\]
since $\lim_{n\rightarrow +\infty}\left( n^{1/2}\,b_n^{-2}\,\tau_n^2\right)=0$, it follows  $\left|U_{n,k,\ell}^{(2)}\right|=o_p(1)$.
In the same way, using in addition  the inequality  $\vert R_{b_n,k}(Y_i)\vert\leq \vert R_{k}(Y_i)\vert$, we get
\[
\left|U_{n,k,\ell}^{(3)}\right| \leq A_2^2 \tau_n^2\,\sqrt{n}\,b_n^{-2}\,\left(A_2^2b_n^{-2}\tau_n^2 +4A_2b_n^{-1}\tau_n+ 4\right) \frac{1}{n}\sum_{i=1}^n\left|R_k(Y_i)R_{\ell}(Y_i)\right|,
\]
\[
\left|U_{n,k,\ell}^{(4)}\right|   \leq A_2^2\left(n^{1/2}b_n^{-2}\,\tau_n^2 \right) \frac{1}{n}\sum_{i=1}^n \left|R_k(Y_i)R_{\ell}(Y_i)\right|, 
\]
and we deduce that $\left|U_{n,k,\ell}^{(3)}\right|=o_p(1)$ and $\left|U_{n,k,\ell}^{(4)}\right|=o_p(1)$.

 \hfill $\square$

\bigskip

\noindent In the following, we define the functions:
\begin{eqnarray*}\label{b70}
A_{k\ell}^{(1)}(y) &=& \frac{g_k(y)g_{\ell}(y)}{f_{b_n}^2(y)}=R_{b_n,k}(y)R_{b_n,\ell}(y),\\ 
A_{k\ell}^{(2)}(y) &=&\frac{g_k(y)\widehat{g}_{\ell , n}(y)+g_{\ell}(y)\widehat{g}_{k , n}(y)}{f_{b_n}^2(y)}=\frac{R_{b_n,k}(y)\widehat{g}_{\ell , n}(y)}{f_{b_n}(y)} + \frac{R_{b_n,\ell}(y)\widehat{g}_{k , n}(y)}{f_{b_n}(y)},\\
A_{k\ell}^{(3)}(y) &=& 2R_{b_n,k}(y)R_{b_n,\ell}(y)\frac{\widehat{f}_{b_n}(y)}{{f}_{b_n}(y)}.
\end{eqnarray*}

\noindent Then, we first have:

\bigskip

\begin{Lemma}\label{a1}
Under the assumptions \ref{ass8} and \ref{ass9}, we have: 
\[
\sqrt{n}\left|\E\left(A_{k\ell}^{(1)}(Y)\right) - \E\left(R_k(Y)R_{\ell}(Y)\right)\right|=o(1).
\]
\end{Lemma}
\noindent\textit{Proof.} The proof is identical to  (4.17) in Zhu and Fang\cite{Zhu1}.
 \hfill $\square$ 

\bigskip

\begin{Lemma}\label{a2} 
Under the assumptions  \ref{ass3}, \ref{ass6}, \ref{ass7},   \ref{ass8} and \ref{ass9}, we have:
\[
\sqrt{n}\,\E\left[A_{k\ell}^{(2)}(Y)\right]=2\sqrt{n}\,\E\left[R_{\ell}(Y)R_{k}(Y)\right]+o(1).
\]
\end{Lemma}
\noindent\textit{Proof.}
With the assumptions \ref{ass8} and \ref{ass9}, we have
\begin{eqnarray*}
\sqrt{n}\,\E\left[\frac{R_{b_n,k}(Y)\widehat{g}_{\ell , n}(Y)}{f_{b_n}(Y)}\right] &=& \sqrt{n}\,\E\left[\frac{R_{b_n,k}(Y)}{f_{b_n}(Y)}2^{j_n}X_{1\ell}K\left(2^{j_n}Y,2^{j_n}Y_1\right)\right].
\end{eqnarray*}
Since $(X_1,Y_1)$ and $Y$ are independent, it follows
\begin{eqnarray*}
& &\sqrt{n}\,\E\left[\frac{R_{b_n,k}(Y)\widehat{g}_{\ell , n}(Y)}{f_{b_n}(Y)}\right]\\
&=& 2^{j_n}\sqrt{n}\int \int\int \frac{R_{b_n,k}(y)}{f_{b_n}(y)} x K\left(2^{j_n}z,2^{j_n}y\right)f_{_{(Y_1,X_{1\ell})}}\left(z,x\right)f(y)\,\,dx\,dz\,dy\\
&=& \sqrt{n}\int\int K\left(2^{j_n}y-u, 2^{j_n}y\right)R_{\ell}\left(y-2^{-j_n}u\right)f\left(y-2^{-j_n}u\right)\frac{R_{b_n,k}(y)f(y)}{f_{b_n}(y)}\,\,du\,dy\\
&=& \mathcal{I}_n \,+\, \mathcal{J}_n, 
\end{eqnarray*}
where
\begin{eqnarray*}
\mathcal{I}_n &=&\sqrt{n}\int\int K\left(2^{j_n}y,2^{j_n}y-u\right)\left(R_{\ell}\left(y-2^{-j_n}u\right)f\left(y-2^{-j_n}u\right)-R_{\ell}(y)f(y)\right)\\
& &\hspace{5cm}\times\, \frac{R_{b_n,k}(y)f(y)}{f_{b_n}(y)}\,\,du\,dy
\end{eqnarray*}
and $\mathcal{J}_n =\sqrt{n}\,\E\left(R_{\ell}(Y)R_{k}(Y)\frac{f^2(Y)}{f_{b_n}^2(Y)}\right)$. 
Furthermore, by Assumption \ref{ass9}, we have
\begin{eqnarray*}
\sqrt{n}\left|\E\left[R_{\ell}(Y)R_{k}(Y)\frac{f^2(Y)}{f_{b_n}^2(Y)}\right] - \E\left[R_{\ell}(Y)R_{k}(Y)\right]\right|  \leqslant \sqrt{n}\left|\E\left[R_{\ell}(Y)R_{k}(Y)\,\mathbf{1}_{\{f(Y)<b_n\}}\right]\right|=o(1).
\end{eqnarray*}
Then, when $n$ is large enough,   $\mathcal{J}_n=\sqrt{n}\,\E\left[R_{\ell}(Y)R_{k}(Y)\right]+o(1)$. For having  the  required result   it is enough   to show that $\mathcal{I}_n = o(1)$. First,
\begin{eqnarray*}
 \left|\mathcal{I}_n\right| &\leqslant & b_n^{-1}\sqrt{n}\left|\int\left(\int K\left(2^{j_n}y,2^{j_n}y+u\right)\left(g_{\ell}\left(y+2^{-j_n}u\right)-g_{\ell}(y)\right)du \right)R_{b_n,k}(y)f(y)dy\right|\\
                            &= & b_n^{-1}\sqrt{n}\left|\int\left(\int K\left(2^{j_n}y,2^{j_n}x\right)\left(g_{\ell}\left(x\right)-g_{\ell}(y)\right)dx \right)R_{b_n,k}(y)f(y)dy\right|.
\end{eqnarray*}
Using  Lemma \ref{lemma41} and the inequality $\left|R_{b_n,k}(y)\right| \leqslant \left|R_k(y)\right|$, we get
$$
\left|\mathcal{I}_n\right|\leqslant  C b_n^{-1}\sqrt{n}2^{-5j_n}\int \left|R_{k}(y)\right|f(y)dy = O\left(b_n^{-1}\sqrt{n}2^{-5j_n}\right).
$$
From Assumption \ref{ass8}, $b_n^{-1}\sqrt{n}2^{-5j_n}\sim n^{c_2+1/2-4c_1} \leqslant n^{-(5c_1+c_2+1/2)}\leqslant n^{-(1/8 + c_2/4)}$;
thus $\mathcal{I}_n = o(1)$.

\hfill $\square$ 

\bigskip

\begin{Lemma}\label{a3}  
Under the assumptions   \ref{ass3}, \ref{ass6}, \ref{ass7},   \ref{ass8} and \ref{ass9}, we have:
\[
\sqrt{n}\,\E\left[A_{k\ell}^{(3)}(Y)\right]=\sqrt{n}\,\E\left[R_{\ell}(Y)R_{k}(Y)\right]+o(1).
\]
\end{Lemma}
\noindent\textit{Proof.} 
Clearly, 
$$\widehat{f}_{b_n}(Y) = V_{n}^{(1)}(Y)+V_{n}^{(2)}(Y)+V_{n}^{(3)}(Y)+V_{n}^{(4)}(Y)+V_{n}^{(5)}(Y),$$
where
\begin{eqnarray*}
 V_{n}^{(1)}(Y) &=& \widehat{f}_n(Y)\mathbf{1}_{\{f(Y)\geqslant  b_n+ C_2\rho_n\}}, \,\,
 V_{n}^{(2)}(Y) = f(Y)\left\{\mathbf{1}_{\{\widehat{f}_n(Y)\geqslant  b_n\}}- \mathbf{1}_{\{f(Y)\geqslant b_n+C_2\rho_n\}}\right\},\\    
 V_{n}^{(3)}(Y) &=& \left(\widehat{f}_n(Y)-f(Y)\right)\left\{\mathbf{1}_{\{\widehat{f}_n(Y)\geqslant b_n\}}- \mathbf{1}_{\{f(Y)\geqslant b_n+C_2\rho_n\}}\right\},\\
 V_{n}^{(4)}(Y) &=& b_n\mathbf{1}_{\{f(Y)< b_n-C_2\rho_n\}},\,\,
 V_{n}^{(5)}(Y) = b_n\left\{\mathbf{1}_{\{\widehat{f}_n(Y)< b_n\}}- \mathbf{1}_{\{f(Y)< b_n-C_2\rho_n\}}\right\},
\end{eqnarray*}
$C_2>0$, and $\rho_n= \left(\log n\right)^{\beta}\times 2^{j_n}n^{-1/2}\,+\,2^{-4j_n}$.
It is then enough to show that
\begin{equation}\label{b200}
\sqrt{n}\E\left[\frac{R_{b_n,k}(Y)R_{b_n,\ell}(Y)V_{n}^{(1)}(Y)}{f_{b_n}(Y)}\right]=\sqrt{n}\E\left[R_k(Y)R_{\ell}(Y)\right]+o(1),
\end{equation}
and that for any   $t\in\{2,\cdots,5\}$,  
\begin{equation}\label{b200bis}
\sqrt{n}\E\left[\frac{R_{b_n,k}(Y)R_{b_n,\ell}(Y)V_{n}^{(t)}(Y)}{f_{b_n}(Y)}\right]\,=\, o(1).
\end{equation}
The proof of  (\ref{b200bis}) is similar to that of (4.23) in Zhu and Fang\cite{Zhu1}. It remains to prove (\ref{b200}). We have:
\begin{eqnarray}  \label{first}
 & & \sqrt{n}\,\E\left[\frac{R_{b_n,k}(Y)R_{b_n,\ell}(Y)V_{n}^{(1)}(Y)}{f_{b_n}(Y)}\right] \nonumber\\
                                                                                         &=& \sqrt{n}\,\E\left[\frac{R_k(Y)R_{\ell}(Y)\widehat{f}_n(Y)}{f(Y)}\mathbf{1}_{\{f(Y)> b_n+ C_2\rho_n\}}\right]\nonumber\\				&=& \sqrt{n}\,\E\left[\frac{R_k(Y)R_{\ell}(Y)}{f(Y)}2^{j_n}K\left(2^{j_n}Y_1,2^{j_n}Y\right)\mathbf{1}_{\{f(Y)> b_n+ C_2\rho_n\}}\right]\nonumber\\	
&= &  \sqrt{n}\int\int_{\{f(y)> b_n+ C_2\rho_n\}}\frac{R_k(y)R_{\ell}(y)}{f(y)}2^{j_n}K\left(2^{j_n}y_1,2^{j_n}y\right)f(y_1)f(y)\,\,dy_1\,dy\nonumber\\
& =&\sqrt{n}\int\int_{\{f(y)> b_n+ C_2\rho_n\}} R_k(y)R_{\ell}(y)K\left(2^{j_n}y+u,2^{j_n}y\right)f(y)du\,dy\nonumber\\
& &+ \,\sqrt{n}\int\int_{\{f(y)> b_n+ C_2\rho_n\}} R_k(y)R_{\ell}(y)K\left(2^{j_n}y+u,2^{j_n}y\right)\left(f(y+u2^{-j_n})-f(y)\right)\,\,du\,dy\nonumber\\
&= &\sqrt{n}\int_{\{f(y)> b_n+ C_2\rho_n\}} R_k(y)R_{\ell}(y)f(y)\,dy \\
& & + \,\sqrt{n}\int\int_{\{f(y)> b_n+ C_2\rho_n\}} R_k(y)R_{\ell}(y)K\left(2^{j_n}y+u,2^{j_n}y\right)\left(f(y+u2^{-j_n})-f(y)\right)du\,dy.\nonumber
\end{eqnarray} 
Further, using Assumption \ref{ass9}-(2) we get
\begin{eqnarray}\label{prem}
\sqrt{n}\int_{\{f(y)> b_n+ C_2\rho_n\}} R_k(y)R_{\ell}(y)f(y)\,dy &= &\sqrt{n}\,\E\left[R_k(Y)R_{\ell}(Y)\right]\nonumber\\
& & -\sqrt{n}\,\E\left[R_k(Y)R_{\ell}(Y)\mathbf{1}_{\{f(Y)\leq  b_n+ C_2\rho_n\}}\right]\nonumber\\
&=&\sqrt{n}\,\E\left[R_k(Y)R_{\ell}(Y)\right]+o(1).
\end{eqnarray}
On the other hand, by a Taylor expansion of  $f$, assumptions \ref{ass3} and \ref{ass7}-(2),  we obtain
\begin{eqnarray*}
& &\sqrt{n}\int\int_{\{f(y)> b_n+ C_2\rho_n\}}  R_k(y)R_{\ell}(y)K\left(2^{j_n}y+u,2^{j_n}y\right)\left(f(y+u2^{-j_n})-f(y)\right)du\,dy\\
& &\hspace{1cm} =\,\sqrt{n}\int\int_{\{f(y)> b_n+ C_2\rho_n\}}  R_k(y)R_{\ell}(y)K\left(2^{j_n}y+u,2^{j_n}y\right)\\
& &\hspace{4cm}\times \left(\frac{2^{-3j_n}}{2}\int_0^1 u^3(1-v)^{2}f^{(3)}\left(y+vu2^{-j_n}\right)dv\right)du\,dy\\
& &\hspace{1cm} =\,\sqrt{n}\frac{2^{-3j_n}}{2}\int\int_{\{f(y)> b_n+ C_2\rho_n\}} R_k(y)R_{\ell}(y)K\left(2^{j_n}y+u,2^{j_n}y\right)\\
& &\hspace{4cm}\times \int_0^1 u^3(1-v)^{2}\left(f^{(3)}\left(y+vu2^{-j_n}\right) -f^{(3)}\left(y\right)\right) dv\,du\,dy.
\end{eqnarray*}
Hence
\begin{eqnarray*}
& &\left|\sqrt{n}\int\int_{\{f(y)> b_n+ C_2\rho_n\}}  R_k(y)R_{\ell}(y)K\left(2^{j_n}y+u,2^{j_n}y\right)\left(f(y+u2^{-j_n})-f(y)\right)du\,dy\right|\\
& & \leqslant  \,\sqrt{n}\frac{c\,2^{-4j_n}}{2}\left(\int\Phi (u)|u|^{4}\,du\right) \int_{\{f(y)> b_n+ C_2\rho_n\}} \frac{\left|R_k(y)R_{\ell}(y)\right|}{f(y)}f(y)\,dy.\\
& &\leqslant \frac{C_3\sqrt{n}\,2^{-4j_n}}{2(b_n+ C_2\rho_n)}\int_{\{f(y)> b_n+ C_2\rho_n\}} \left|R_k(y)R_{\ell}(y)\right|f(y)\,dy\\
& &\leqslant \frac{C_3b_n^{-1}\sqrt{n}\,2^{-4j_n}}{2(1+C_2\, b_n^{-1}\rho_n)}\E\left(  \left|R_k(Y)R_{\ell}(Y)\right|\right).
\end{eqnarray*}
However, $b_n^{-1}\sqrt{n}\,2^{-4j_n}\sim n^{-4(c_1-c_2/4-1/8)}$ and since  $c_1-c_2/4-1/8>0$ (see Assumption  \ref{ass8}), it follows: $\lim_{n\rightarrow +\infty}\left(b_n^{-1}\sqrt{n}\,2^{-4j_n}\right)=0$. Furthermore, $b_n^{-1}\rho_n=b_n^{-1}2^{-4j_n}\left(1+\varepsilon_n\right)$, where $\varepsilon_n=n^{-1/2}2^{5j_n}\left(\log n\right)^\beta\sim n^{5(c_1-1/10)}\left(\log n\right)^\beta\rightarrow 0$ as $n\rightarrow +\infty$. Consequently,  $b_n^{-1}\rho_n\sim b_n^{-1}2^{-4j_n}\sim n^{4(c_2/4-c_1)}$. Since $c_2/4<1/8+c_2/4<c_1$ by Assumption \ref{ass8}, we deduce that $\lim_{n\rightarrow +\infty}\left(b_n^{-1}\rho_n\right)=0$.
Finally, we get:
\begin{equation*}\label{b202}
\sqrt{n}\int\int_{\{f(Y)> b_n+ C_2\rho_n\}} R_k(y)R_{\ell}(y)K\left(2^{j_n}y,2^{j_n}y+u\right)\left(f(y+u2^{-j_n})-f(y)\right)du\,dy=o(1);
\end{equation*}
this equality and  the equations (\ref{first}),   (\ref{prem}) allow to conclude that (\ref{b200}) holds.
\hfill $\square$ 
\bigskip

\noindent In all what follows, we consider for $j\in\{1,2,3\}$:
\[
\mathcal{E}_{k\ell}^{(j)}=\frac{1}{\sqrt{n}}\sum_{i=1}^n\Big\{A_{k\ell}^{(j)}\left(Y_i\right)-\E\left( A_{k\ell}^{(j)}\left(Y\right)\right)\Big\}.
\]
We first have:

\bigskip

\begin{Lemma}\label{b31}
Under assumptions \ref{ass1}, \ref{ass3}, \ref{ass6}, \ref{ass7}, \ref{ass8} and \ref{ass9}, we have:    
\[
\mathcal{E}_{k\ell}^{(1)} =\frac{1}{\sqrt{n}}\sum_{i=1}^n\Big\{R_{k}\left(Y_i\right)R_{\ell}\left(Y_i\right)-\E\left(R_{k}\left(Y\right)R_{\ell}\left(Y\right)\right)\Big\} + o_p(1).
\] 
\end{Lemma}
\noindent\textit{Proof.} The proof is identical to that of step 2 of  the proof of  Theorem 2.1 in Zhu and Fang\cite{Zhu1}.
\hfill $\square$ 

\bigskip

\begin{Lemma}\label{b26}
Under assumptions \ref{ass1}, \ref{ass3}, \ref{ass4}, \ref{ass6}, \ref{ass7}, \ref{ass8} and \ref{ass9}, we have    
\begin{multline*}
\mathcal{E}_{k\ell}^{(2)}
=\frac{1}{\sqrt{n}}\sum_{i=1}^n\left\{R_{b_n,k}\left(Y_i\right)R_{b_n,\ell}\left(Y_i\right)+ \frac{1}{2}X_{\ell i}R_{b_n,k}(Y_i)\frac{f(Y_i)}{f_{b_n}(Y_i)}+\frac{1}{2}X_{ki}R_{b_n,\ell}(Y_i)\frac{f(Y_i)}{f_{b_n}(Y_i)}\right\}\\
 -\,\E\left\{ R_{b_n,k}\left(Y\right)R_{b_n,\ell}\left(Y\right)+ \frac{1}{2}X_{\ell}R_{b_n,k}(Y)\frac{f(Y)}{f_{b_n}(Y)}+\frac{1}{2}X_{k}R_{b_n,\ell}(Y)\frac{f(Y)}{f_{b_n}(Y)} \right\} + o_p(1)
\end{multline*}    
\end{Lemma}
\noindent\textit{Proof.}
Applying the similar argument used in several works in order to approximate a sum by a U-statistic (\cite{Hardle1}, \cite{Hardle2}, \cite{Nolan}, \cite{Powell1}, \cite{Stone}), we obtain
\begin{multline}\label{b25}
\frac{1}{\sqrt{n}}\sum_{i=1}^n\left\{R_{b_n,k}\left(Y_i\right)\frac{\widehat{g}_{\ell , n}\left(Y_i\right)}{f_{b_n}\left(Y_i\right)}-\E\left[R_{b_n,k}\left(Y\right)\frac{\widehat{g}_{\ell , n}\left(Y\right)}{f_{b_n}\left(Y\right)}\right]\right\}\\
= \frac{1}{\sqrt{n}}\sum_{i=1}^n\mathcal{V}_{n,k,i}-\E\left(\mathcal{V}_{n,k,i}\right)+o_p(1),
\end{multline}
where
\[
\mathcal{V}_{n,k,i}=\frac{1}{2}\int\int 2^{j_n}K\left(2^{j_n}y,2^{j_n}Y_i\right)\bigg(\frac{x\,R_{b_n,k}(Y_i)}{f_{b_n}(Y_i)}+\frac{X_{\ell i}R_{b_n,k}(y)}{f_{b_n}(y)}\bigg) f_{_{(X_{\ell},Y)}}(x,y)dx\,dy.
\]      
From the equalities
\[
\int x \,\frac{f_{(X_{\ell},Y)}(x,y)}{f(y)}\,dx= R_\ell(y),\,\,\,
\int  f_{(X_{\ell},Y)}(x,y )\,dx=f(y),\,\,\, R_{b_n,k}(Y_i)=R_k(Y_i)\,\varepsilon_{b_n}(Y_i),
\]
where $\varepsilon_{b_n}(y)= f(y)/f_{b_n}(y)$,  and the property $\int K(x,y)\,dy=1$,        it follows    
\begin{eqnarray*}                   
\mathcal{V}_{n,k,i}& =&  \frac{R_{b_n,k}(Y_i)}{2f_{b_n}(Y_i)}\int 2^{j_n}K\left(2^{j_n}y,2^{j_n}Y_i\right)R_{\ell}(y)f(y)dy\\
&&   +\frac{X_{\ell i}}{2}\int 2^{j_n}K\left(2^{j_n}y,2^{j_n}Y_i\right) \frac{R_{b_n,k}(y)}{f_{b_n}(y)}f(y)dy\\
&=&\frac{1}{2}R_k(Y_i)R_{\ell}(Y_i)\varepsilon_{_{b_n}}^2(Y_i)+\frac{1}{2}X_{\ell i}R_k(Y_i)\varepsilon_{_{b_n}}^2(Y_i)\,+\,U_n^{(5)}(Y_i)\,+\,U_n^{(6)}(Y_i),
\end{eqnarray*}
where 
\[
U_n^{(5)}(Y_i) = \frac{R_{b_n,k}(Y_i)}{2f_{b_n}(Y_i)}\int 2^{j_n}K\left(2^{j_n}y,2^{j_n}Y_i\right)\Big(R_{\ell}(y)f(y)-R_{\ell}(Y_i)f(Y_i)\Big)\,\,dy,
\]
\[
U_n^{(6)}(Y_i)  =\frac{X_{\ell i}}{2}\int 2^{j_n}K\left(2^{j_n}y,2^{j_n}Y_i\right)\left( \frac{R_{b_n,k}(y)f(y)}{f_{b_n}(y)}-\frac{R_{b_n,k}(Y_i)f(Y_i)}{f_{b_n}(Y_i)}\right)\,\,dy.
\]
However, from  Lemma \ref{lemma41}, 
\begin{eqnarray*}
& & \E\bigg(\left(U_n^{(5)}(Y)\right)^2\bigg)\\
                      & = & \left(\frac{2^{j_n}}{2}\right)^2\int \left(\frac{R_{b_n,k}(y)}{f_{b_n}(y)}\right)^2\left(\int K\left(2^{j_n}x,2^{j_n}y\right)\Big(g_{\ell}(x)-g_{\ell}(y)\Big)\,\,dx\right)^2f(y)\,\,dy \\
                     & \leqslant & \frac{C^2_1}{4}\,2^{-8j_n}b_n^{-2}\int R_k(y)^2f(y)dy=\frac{C^2_1}{4}\,2^{-8j_n}\,b_n^{-2}\E\left(R_k(Y)^2\right).
\end{eqnarray*}
From Assumption \ref{ass9}-(1), $\E\left(R_k(Y)^2\right)<+\infty$, and   from Assumption \ref{ass8},  $\quad 2^{-8j_n}\,b_n^{-2} \sim n^{-2(4c_1-c_2)}$. Since $4c_1-c_2>0$ because $c_1>1/8+c_2/4>c_2/4$, it follows that $\E\left(\left(U_n^{(5)}(Y)\right)^2\right)=o(1)$.                  
Since $Var\left(n^{-1/2}\sum_{i=1}^nU^{(5)}_n(Y_i)\right) =Var\left( U^{(5)}_n(Y)\right)\leq \E\bigg(\left(U_n^{(5)}(Y)\right)^2\bigg) $, it follows that $n^{-1/2}\sum_{i=1}^n\left(U^{(5)}_n(Y_i)-\E (U^{(5)}_n(Y_i))\right) =o_p(1)$. On the other hand,
\begin{eqnarray*}
& &\E\left[\left(U_n^{(6)}(Y)\right)^2\right]\\
 &=& \E\left[\left( \frac{X_{\ell}}{2}\int 2^{j_n}K\left(2^{j_n}x,2^{j_n}Y\right)\left(\frac{R_{b_n,k}(x)f(x)}{f_{b_n}(x)}-\frac{R_{b_n,k}(Y)f(Y)}{f_{b_n}(Y)}\right)\,\,dx\right)^2\right] \\
                                           &=& \frac{1}{4}\,\int\int \left(t \int 2^{j_n}K\left(2^{j_n}x,2^{j_n}y\right)\left(\frac{R_{b_n,k}(x)f(x)}{f_{b_n}(x)}-\frac{R_{b_n,k}(y)f(y)}{f_{b_n}(y)}\right)\,\,dx\right)^2 f_{(X_\ell,Y)} ( t,y)\,dt\,dy\\
 &=& \frac{1}{4}\,\int\left(\int t^2 \frac{f_{(X_\ell,Y)} ( t,y)}{f(y)} dt\right)\\
&&\times \left( \int 2^{j_n}K\left(2^{j_n}x,2^{j_n}y\right)\left(\frac{R_{b_n,k}(x)f(x)}{f_{b_n}(x)}-\frac{R_{b_n,k}(y)f(y)}{f_{b_n}(y)}\right)\,\,dx\right)^2f(y)\,\,dy \\
&=& \frac{1}{4}\,\E \left[ \E\left(X_{\ell}^2|Y\right)\,\left( \int 2^{j_n}K\left(2^{j_n}x,2^{j_n}Y\right)\left\{ \frac{R_{b_n,k}(x)f(x)}{f_{b_n}(x)}-\frac{R_{b_n,k}(Y)f(Y)}{f_{b_n}(Y)}\right\}dx\right)^2\right]\\
                                           &=& \frac{1}{4}\,\E \left[\E\left(X_{\ell}^2|Y\right)\,\left(\int K\left(2^{j_n}Y,2^{j_n}Y+u\right)\right.\right.\\
&&\left.\left.\times\left( \frac{R_{b_n,k}(Y+2^{-j_n}u)f(Y+2^{-j_n}u)}{f_{b_n}(Y+2^{-j_n}u)}-\frac{R_{b_n,k}(Y)f(Y)}{f_{b_n}(Y)}\right)du\right)^2\right]\\
&=& \frac{1}{4}\,\E \left[\E\left(X_{\ell}^2|Y\right)\,h_n^2(Y)\right],
\end{eqnarray*}
where $h_n(Y)=\int w_n(Y,u)\,du$ with $w_n(Y,u)=  K\left(2^{j_n}Y,2^{j_n}Y+u\right)\,m_n(Y,u)$ and 
\[
m_n(Y,u)=   \frac{R_{b_n,k}(Y+2^{-j_n}u)f(Y+2^{-j_n}u)}{f_{b_n}(Y+2^{-j_n}u)}-\frac{R_{b_n,k}(Y)f(Y)}{f_{b_n}(Y)}.
\]
From continuity of $f$ and $g_k$ and the fact that the sequence  $\left(K\left(2^{j_n}Y,2^{j_n}Y+u\right)\right)_{n\in\mathbb{N}}$ is bounded (see Assumption \ref{ass7}-(1)),  we obtain $\lim_{n\rightarrow +\infty}\left(w_n(Y,u)\right)=0$. Further, since  $\varepsilon_{b_n}=f/f_{b_n}\leq 1$ and $R_{b_n,k}=R_k\,\varepsilon_{b_n}$,  it follows 
\begin{eqnarray*}
\vert m_n(Y,u)\vert&\leq& \bigg\vert \bigg(R_{k}(Y+2^{-j_n}u)-R_k(Y)\bigg)\,\varepsilon_{b_n}^2(Y+2^{-j_n}u)\\
& &+ R_k(Y) \bigg(\varepsilon_{b_n}^2(Y+2^{-j_n}u)-\varepsilon_{b_n}^2(Y)\bigg)\bigg\vert\\
&\leq&  \bigg\vert R_{k}(Y+2^{-j_n}u)-R_k(Y)\bigg\vert+2\vert  R_k(Y)  \vert\\
&\leq& c2^{-j_n}\vert u\vert +2\vert  R_k(Y)  \vert.
\end{eqnarray*}
The sequence $\left( 2^{-j_n}\right)_{n\in\mathbb{N}}$ being   bounded, there exists  $M >0$   such that $\vert m_n(Y,u)\vert\leq M\vert u\vert +2\vert  R_k(Y)  \vert$.  Hence
\[
\int \vert w_n(Y,u)\vert\,du\leq M\int\vert u\vert\Phi(u)\,du +2\vert  R_k(Y)  \vert\int\Phi(u)\,du<+\infty;
\]
and using the dominated convergence theorem we get: $\lim_{n\rightarrow +\infty}\left(h_n(Y)\right)=0$. On the other hand,
\begin{eqnarray*}
h_n^2(Y)&\leq&\int  w_n^2(Y,u)\,du=\int  K^2\left(2^{j_n}Y,2^{j_n}Y+u\right)\,m_n^2(Y,u)\,du\\
&\leq&\int  K^2\left(2^{j_n}Y,2^{j_n}Y+u\right)\,\left(M\vert u\vert +2\vert  R_k(Y)  \vert\right)^2\,du\\
&\leq&2M^2\int  K^2\left(2^{j_n}Y,2^{j_n}Y+u\right)\,u^2\,du+8 R^2_k(Y)\int  K^2\left(2^{j_n}Y,2^{j_n}Y+u\right) \,du.
\end{eqnarray*}
Then, using Assumption \ref{ass7}-(1) and the property $\int K^2(x,y)\,dy\leqslant D^2$ (see Remark  \ref{rmq2}), we get $h_n^2(Y)\leq M_1+M_2R^2_k(Y)$, where $M_1=2M^2\int u^2\Phi^2(u)\,du$ and $M_2=8D^2$.
Therefore,  $\vert E\left(X_{\ell}^2|Y\right)\,h_n^2(Y)\vert\leq \phi (Y)$, where 
$
\phi (Y)=M_1\E\left(X_{\ell}^2|Y\right)+ M^2_2\E\left(X_{\ell}^2R_k^2(Y)|Y\right)$.
This random variable is integrable since
\begin{eqnarray*}
\mathbb{E}\left(\vert\phi (Y)\vert\right)&=&M_1\E\left(X_{\ell}^2 \right)+ M_2\E\left(X_{\ell}^2 R_k^2(Y) \right)\\
&\leqslant&M_1G^2+M_2G^2 \E\left( R_k^2(Y) \right) <+\infty.
\end{eqnarray*}
Applying the dominated convergence theorem, we obtain
\[
\lim_{n\rightarrow +\infty}\E\left[\left(U_n^{(6)}(Y)\right)^2\right]=\frac{1}{4}\,\E \left[\E\left(X_{\ell}^2|Y\right)\,\lim_{n\rightarrow +\infty}\left(h_n^2(Y)\right)\right]=0
\]
and we deduce, since $Var\left(n^{-1/2}\sum_{i=1}^nU^{(6)}_n(Y_i)\right) =Var\left( U^{(6)}_n(Y)\right)\leq \E\bigg(\left(U_n^{(6)}(Y)\right)^2\bigg) $,  that $n^{-1/2}\sum_{i=1}^n\left(U^{(6)}_n(Y_i)-\E (U^{(6)}_n(Y_i))\right) =o_p(1)$.  These results allow to conclude that
\begin{eqnarray*}
 & &\frac{1}{\sqrt{n}}\sum_{i=1}^n \mathcal{V}_{n,k,i} - \E\left(\mathcal{V}_{n,k,i}\right)\\
	    & &\hspace{2cm}  =\,\frac{1}{\sqrt{n}}\sum_{i=1}^n \bigg\{\frac{1}{2}R_k\left(Y_i\right)R_{\ell}\left(Y_i\right)\varepsilon_{b_n}^2\left(Y_i\right)+\frac{1}{2}X_{\ell i}R_k(Y_i)\varepsilon_{b_n}^2\left(Y_i\right)\\
   & & \hspace{2cm} -\,\E\left(\frac{1}{2}R_k\left(Y\right)R_{\ell}\left(Y\right)\varepsilon_{b_n}^2\left(Y\right)+\frac{1}{2}X_{\ell i}R_k(Y)\varepsilon_{b_n}^2\left(Y\right) \right)\bigg\}+o_p(1),
\end{eqnarray*}
and from Eq.(\ref{b25}), we obtain
\begin{multline}\label{b119} 
\frac{1}{\sqrt{n}}\sum_{i=1}^n\bigg\{R_{b_n,k}\left(Y_i\right)\frac{\widehat{g}_{\ell , n}\left(Y_i\right)}{f_{b_n}\left(Y_i\right)}-\E\left[R_{b_n,k}\left(Y\right)\frac{\widehat{g}_{\ell , n}\left(Y\right)}{f_{b_n}\left(Y\right)}\right]\bigg\}\\
=\frac{1}{\sqrt{n}}\sum_{i=1}^n \bigg\{\frac{1}{2}R_k\left(Y_i\right)R_{\ell}\left(Y_i\right)\varepsilon_{b_n}^2\left(Y_i\right)+\frac{1}{2}X_{\ell i}R_k(Y_i)\varepsilon_{b_n}^2\left(Y_i\right)\\ -\E\left(\frac{1}{2}R_k\left(Y\right)R_{\ell}\left(Y\right)\varepsilon_{b_n}^2\left(Y\right)+\frac{1}{2}X_{\ell i}R_k(Y)\varepsilon_{b_n}^2\left(Y\right) \right)\bigg\}+o_p(1)\\
=\frac{1}{\sqrt{n}}\sum_{i=1}^n \bigg\{\frac{1}{2}R_{b_n,k}\left(Y_i\right)R_{b_n,\ell}\left(Y_i\right)+\frac{1}{2}X_{\ell i}R_{b_n,k}(Y_i)\frac{f(Y_i)}{f_{b_n}(Y_i)}\\ -\E\left(\frac{1}{2}R_{b_n,k}\left(Y\right)R_{b_n,\ell}\left(Y\right)+\frac{1}{2}X_{\ell i}R_{b_n,k}(Y)\frac{f(Y)}{f_{b_n}(Y)} \right)\bigg\}+o_p(1).
\end{multline}
By replacing  $k$  with $\ell$  and adding the results, we obtain the   required result.
 
\hfill $\square$

\begin{Lemma}\label{b33}
Under assumptions \ref{ass1}, \ref{ass3},  \ref{ass6}, \ref{ass7}, \ref{ass8} and \ref{ass9}, we have    
\begin{multline*}\label{b27}
\mathcal{E}_{k\ell}^{(3)}
=\frac{2}{\sqrt{n}}\sum_{i=1}^n\bigg\{R_{b_n,k}\left(Y_i\right)R_{b_n,\ell}\left(Y_i\right)\frac{f\left(Y_i\right)}{f_{b_n}\left(Y_i\right)}-\E\left[R_{b_n,k}\left(Y\right)R_{b_n,\ell}\left(Y\right)\frac{f\left(Y_i\right)}{f_{b_n}\left(Y\right)}\right]\bigg\} + o_p(1).
\end{multline*}
 \end{Lemma}
\noindent\textit{Proof.} The proof  similar  to that of step 2 of  the proof of  Theorem 2.1 in Zhu and Fang\cite{Zhu1}.
\hfill $\square$

\begin{Lemma}\label{b34}
Under assumptions \ref{ass1}, \ref{ass3},  \ref{ass6}, \ref{ass7}, \ref{ass8} and \ref{ass9}, we have
 
$
\frac{1}{\sqrt{n}}\sum_{i=1}^n\Bigg\{R_{b_n,k}\left(Y_i\right)R_{b_n,\ell}\left(Y_i\right)+\frac{X_{\ell i}R_{b_n,k}\left(Y_i\right)f(Y_i)}{2f_{b_n}(Y_i)}
+\frac{X_{ki}R_{b_n,\ell}\left(Y_i\right)f(Y_i)}{2f_{b_n}(Y_i)}$

$\hspace{1cm} -\,\E\left[R_{b_n,k}\left(Y\right)R_{b_n,\ell}\left(Y\right)+\frac{X_{\ell}R_{b_n,k}\left(Y\right)f(Y)}{2f_{b_n}(Y)}
+\frac{X_{k}R_{b_n,\ell}\left(Y\right)f(Y)}{2f_{b_n}(Y)}\right]
\Bigg\}$

$\hspace{2cm}  =\frac{1}{\sqrt{n}}\sum_{i=1}^n\bigg\{R_{k}\left(Y_i\right)R_{\ell}\left(Y_i\right)+\frac{1}{2}X_{\ell i}R_{k}\left(Y_i\right)+\frac{1}{2}X_{ki}R_{\ell}\left(Y_i\right)
-2\E\left[R_{k}\left(Y\right)R_{\ell}\left(Y\right)\right]\bigg\}+o_p(1)
$
\end{Lemma}
\noindent\textit{Proof.} Putting 
\[
\mathcal{Z}_{k,\ell, n}=\frac{1}{\sqrt{n}}\sum_{i=1}^n\left\{\frac{X_{\ell i}R_{b_n,k}\left(Y_i\right)f\left(Y_i\right)}{2f_{b_n}\left(Y_i\right)} - \frac{1}{2}X_{\ell i}R_{k}\left(Y_i\right)  -\E\left[\frac{X_{\ell}R_{b_n,k}\left(Y\right)f\left(Y\right)}{2f_{b_n}\left(Y\right)} - \frac{1}{2}X_{\ell}R_{k}\left(Y\right)\right]\right\},
\]
and using Assumption \ref{ass1}, we get
\begin{eqnarray*}
\E\left(\mathcal{Z}_{k,\ell, n}^2\right) &\leqslant & \frac{G^2}{4}\E \left(\frac{R_{b_n,k}\left(Y\right)f\left(Y\right)}{f_{b_n}\left(Y\right)} - R_{k}\left(Y\right)\right)^2\\
										&\leqslant & \frac{G^2}{4}\int_{\{f(y)<b_n\}}\left(\left|\frac{R_{b_n,k}\left(y\right)f\left(y\right)}{f_{b_n}\left(y\right)}\right| + \left|R_{k}\left(y\right)\right|\right)^2f(y)\,dy\\
										&\leqslant & \frac{G^2}{4}\int_{\{f(y)<b_n\}}\left(\left|R_{b_n,k}\left(y\right)\right| + \left|R_{k}\left(y\right)\right|\right)^2f(y)\,dy\\
										&\leqslant & G^2\E\left[\left(R_{k}\left(Y\right)\right)^2\mathbf{1}_{\{f(y)<b_n\}}\right];\\
\end{eqnarray*}
from Assumption \ref{ass9} we deduce that $\,\E\left(\mathcal{Z}_{k,\ell, n}^2\right) \rightarrow 0 \,$ as $\, n \rightarrow +\infty$. Thus
\begin{equation}\label{b66}
\frac{1}{\sqrt{n}}\sum_{i=1}^n\left\{\frac{X_{\ell i}R_{b_n,k}\left(Y_i\right)f\left(Y_i\right)}{2f_{b_n}\left(Y_i\right)} - \frac{1}{2}X_{\ell i}R_{k}\left(Y_i\right)  -\E\left[\frac{X_{\ell}R_{b_n,k}\left(Y\right)f\left(Y\right)}{2f_{b_n}\left(Y\right)} - \frac{1}{2}X_{\ell}R_{k}\left(Y\right)\right]\right\}\,=\,o_p(1).
\end{equation}
By  inverting $k$   and  $\ell$,   and adding the result to the previous one, we get
\begin{align*}
&\phantom{=} \frac{1}{\sqrt{n}}\sum_{i=1}^n\left\{\frac{X_{ki}R_{b_n,\ell}\left(Y_i\right)f\left(Y_i\right)}{2f_{b_n}\left(Y_i\right)} + \frac{X_{\ell i}R_{b_n,k}\left(Y_i\right)f\left(Y_i\right)}{2f_{b_n}\left(Y_i\right)}\right. \\
&\phantom{=}\left.
\hspace{2cm} - \,\E\left[\frac{R_{\ell}\left(Y\right)R_{b_n,k}\left(Y\right)f\left(Y\right)}{2f_{b_n}\left(Y\right)} + \frac{R_k\left(Y\right)R_{b_n,\ell}\left(Y\right)f\left(Y\right)}{2f_{b_n}\left(Y\right)}\right]\right\} \\ 
&= 
\frac{1}{\sqrt{n}}\sum_{i=1}^n\left\{\frac{1}{2}X_{\ell i}R_{k}\left(Y_i\right) + \frac{1}{2}X_{ki}R_{\ell}\left(Y_i\right) -  \E\left[R_{k}\left(Y\right)R_{\ell}\left(Y\right)\right]\right\}+o_p(1),
\end{align*}
and adding this  result to the one of   Lemma \ref{b31} yields the required result.
\hfill $\square$ 

\bigskip

\begin{Lemma}\label{b68} 
Under assumptions \ref{ass1}, \ref{ass3},  \ref{ass6}, \ref{ass7}, \ref{ass8} and \ref{ass9}, we have
\begin{multline*} \frac{1}{\sqrt{n}}\sum_{i=1}^n\bigg\{R_{b_n,k}\left(Y_i\right)R_{b_n,\ell}\left(Y_i\right)\frac{f(Y_i)}{f_{b_n}(Y_i)}-\E\left[R_{b_n,k}\left(Y\right)R_{b_n,\ell}\left(Y\right)\frac{f(Y)}{f_{b_n}(Y)}\right]\bigg\}\\
=\frac{1}{\sqrt{n}}\sum_{i=1}^n\bigg\{R_k\left(Y_i\right)R_{\ell}\left(Y_i\right)-\E\left[R_k\left(Y\right)R_{\ell}\left(Y\right)\right]\bigg\}+o_p(1).
\end{multline*}
\end{Lemma}
\noindent\textit{Proof.} The proof is obtained by using similar arguments than in the proof of Lemma \ref{b34} from
\begin{multline*} 
\mathcal{Y}_{k,\ell,n}=\frac{1}{\sqrt{n}}\sum_{i=1}^n\bigg\{R_{b_n,k}\left(Y_i\right)R_{b_n,\ell}\left(Y_i\right)\frac{f(Y_i)}{f_{b_n}(Y_i)}-R_k\left(Y_i\right)R_{\ell}\left(Y_i\right)\\
-\E\left[R_{b_n,k}\left(Y\right)R_{b_n,\ell}\left(Y\right)\frac{f(Y)}{f_{b_n}(Y)} - R_k\left(Y\right)R_{\ell}\left(Y\right)\right]\bigg\}.
\end{multline*}
\hfill $\square$ 
 
\subsection{Proof of Theorem \ref{theo1}}
Let us denote by  $\widehat{\lambda}_{k,\ell}^{(n)}$ the $(k,\ell)$-th entry  of the  $d\times d$  matrix $\widehat{\Lambda}_n$.  It is easily seen that
\begin{eqnarray*}
\sqrt{n}\,\widehat{\lambda}_{k,\ell}^{(n)}=\frac{1}{\sqrt{n}}\sum_{i=1}^n\frac{\widehat{g}_{k , n}(Y_i)\widehat{g}_{\ell , n}(Y_i)}{\widehat{f}_{b_n}^2(Y_i)}&=&\frac{1}{\sqrt{n}}\sum_{i=1}^n\Big\{A_{k\ell}^{(1)}(Y_i)+A_{k\ell}^{(2)}(Y_i)-A_{k\ell}^{(3)}(Y_i)\Big\} \\
&&- U_{n,k,\ell}^{(1)}+U_{n,k,\ell}^{(2)} + U_{n,k,\ell}^{(3)} - U_{n,k,\ell}^{(4)},
\end{eqnarray*}
and  from Lemma  \ref{cl14}  we get
\begin{equation}\label{b22} 
  \sqrt{n}\,\widehat{\lambda}_{k,\ell}^{(n)} =\frac{1}{\sqrt{n}}\sum_{i=1}^n\Big\{A_{k\ell}^{(1)}(Y_i)+A_{k\ell}^{(2)}(Y_i)-A_{k\ell}^{(3)}(Y_i)\Big\}+o_p(1).
\end{equation}
Therefore, putting $\nu_{k,\ell}=\mathbb{E}\left(A_{k\ell}^{(1)}(Y)+A_{k\ell}^{(2)}(Y)-A_{k\ell}^{(3)}(Y)\right)$, we have $\sqrt{n}\left(\widehat{\lambda}_{k,\ell}^{(n)}-\nu_{k,\ell}\right)=
\mathcal{E}_{k\ell}^{(1)}+\mathcal{E}_{k\ell}^{(2)}-\mathcal{E}_{k\ell}^{(3)}+\,o_p(1)$. 
Then, from   Lemmas \ref{b31}, \ref{b26}, \ref{b34}, \ref{b33} and \ref{b68}, it follows
\[
\sqrt{n}\left(\widehat{\lambda}_{k,\ell}^{(n)}-\nu_{k,\ell}\right)=\frac{1}{\sqrt{n}}\sum_{i=1}^n\bigg\{\frac{1}{2}\left(X_{\ell i}R_k(Y_i)+X_{ki}R_{\ell}(Y_i)\right)-\E\left(R_k(Y)R_{\ell}(Y)\right)\bigg\}+o_p(1),
\]
and using  Lemmas \ref{a1},  \ref{a2} and \ref{a3}, we get $\sqrt{n}\,\nu_{k,\ell}=\sqrt{n}\,\lambda_{k,\ell}+o(1)$, where $\lambda_{k,\ell}=\E\left(R_k(Y)R_{\ell}(Y)\right)$, that is the $(k,\ell)$-th entry  of the $d\times d$ matrix $\Lambda$. Hence,
\[
\sqrt{n}\left(\widehat{\lambda}_{k,\ell}^{(n)}-\lambda_{k,\ell}\right)=\frac{1}{\sqrt{n}}\sum_{i=1}^n\bigg\{\frac{1}{2}\left(X_{\ell i}R_k(Y_i)+X_{ki}R_{\ell}(Y_i)\right)-\E\left(R_k(Y)R_{\ell}(Y)\right)\bigg\}+o_p(1),
\]
Clearly,
\[
\E\left(\frac{1}{2}\left(X_{\ell }R_k(Y)+X_{k}R_{\ell}(Y)\right)\right)=\E\Big(X_{\ell}R_k(Y)\Big) =\E\Big(R_{\ell}(Y)R_k(Y)\Big),
\]
and, putting  $\mathcal{H}_{n}=\sqrt{n}\,\left(\widehat{\Lambda}_n-\Lambda\right)$ and $\mathcal{H}^{(n)}_{k\ell}=\sqrt{n}\left(\widehat{\lambda}_{k,\ell}^{(n)}-\lambda_{k,\ell}\right)$, we have
\[
tr\left(A^T\mathcal{H}_{n}\right)=\sum_{k=1}^d\sum_{\ell=1}^da_{k\ell}\,\mathcal{H}^{(n)}_{k\ell}
=\frac{1}{\sqrt{n}}\sum_{i=1}^n\left(\mathcal{U}_i-\E\left(\mathcal{U}_i\right)\right)+o_p(1),
\]
where 
\[
\mathcal{U}_i=\sum_{k=1}^d\sum_{\ell=1}^d\frac{a_{k\ell}}{2}\left(X_{\ell i}R_k(Y_i)+X_{ki}R_{\ell}(Y_i)\right).
\]
From the central limit theorem and  Slutsky's theorem we deduce  that $tr\left(A^T\mathcal{H}_{n}\right)\stackrel{\mathscr{D}}{\rightarrow}\mathcal{N}\left(0,\sigma_{A}^{2}\right)$,   as $n\rightarrow +\infty$, where $\sigma_A^2$ is given in (\ref{siga}).    Then, using Levy's theorem, we conclude that  $ \mathcal{H}_{n}\stackrel{\mathscr{D}}{\rightarrow}\mathcal{H}$, as $n\rightarrow +\infty$, where $\mathcal{H}$ has a normal distribution in $\mathscr{M}_d(\mathbb{R})$ with   $tr\left(A^T\mathcal{H}\right)\leadsto  \mathcal{N}\left(0,\sigma_{A}^{2}\right)$.
\hfill $\square$ 

\subsection{Proof of Theorem \ref{theo2}}
It is known from Theorem 2.2 in Zhu and Fang\cite{Zhu1} that  $
\sqrt{n}\left(\widehat{\beta}_j-\beta_j\right)\stackrel{\mathscr{D}}{\rightarrow}\mathcal{G}_j=\sum_{\stackrel{r=1}{r\neq j}}^d(\lambda_j-\lambda_r)^{-1}\beta_r\beta_j^T\mathcal{H}\beta_j$. Since $\mathcal{G}_j$ is a linear function of $\mathcal{H}$, it has a normal distribution $\mathcal{N}(0,\Sigma_j)$, where $\Sigma_j$ is a covariance matrix that will now be specified. Clearly, for any $u=\left(u_1,\cdots,u_d\right)\in\mathbb{R}^d$, we have $u^T\mathcal{G}_j=tr(A\mathcal{H})$, where $A=\beta_ju^T\xi_j\beta_j^T$ with $\xi_j=\sum_{\stackrel{r=1}{r\neq j}}^d(\lambda_j-\lambda_r)^{-1}\,\beta_r$. Then, from  Theorem \ref{theo1} we get $u^T\mathcal{G}_j\leadsto  \mathcal{N}\left(0,\sigma_{A}^{2}\right)$, where $\sigma_A^2=Var\left(\sum_{p=1}^d\sum_{q=1}^d2^{-1}a_{pq}\left(X_{q}R_p(Y)+X_{p}R_{q}(Y)\right)\right)$, and $a_{pq}$ is the $(p,q)$-th entry of $A$. However, $a_{pq}=\sum_{k=1}^d\beta_{jp}u_k\xi_{jk}\beta_{jq}$, and, therefore,
\begin{eqnarray*}
\sum_{p=1}^d\sum_{q=1}^d\frac{a_{pq}}{2}\left(X_{q}R_p(Y)+X_{p}R_{q}(Y)\right)&=&\sum_{k=1}^du_k\xi_{jk}\sum_{p=1}^d\sum_{q=1}^d\frac{\beta_{jp}\beta_{jq}}{2}\left(X_{q}R_p(Y)+X_{p}R_{q}(Y)\right).
\end{eqnarray*}
Since $\xi_{jk}=\sum_{\stackrel{r=1}{r\neq j}}^d(\lambda_j-\lambda_r)^{-1}\,\beta_{rk}$, it follows $\sigma_A^2=Var\left(\sum_{k=1}^du_k\mathcal{W}_{jk}\right)=u^T\Theta_ju$, where $\Theta_j$ is the covariance matrix of $\mathcal{W}_j=\left(\mathcal{W}_{j1},\cdots,\mathcal{W}_{jd}\right)^T$. Since this later equality holds for all $u$ in $\mathbb{R}^d$, we deduce that $\Sigma_j=\Theta_j$.
\hfill $\square$

\end{document}